\documentclass[11pt]{amsart}

\usepackage{soul}

\usepackage{amsmath,amsthm, amssymb, color}
\usepackage{graphicx}
\usepackage{geometry}
\usepackage[normalem]{ulem}

\usepackage{amsmath, amsthm, amssymb, bm, mathtools}
\usepackage{enumerate}

\newcommand{\R}{{\mathbb {R}}}
\newcommand{\Pb}{{\mathbb {P}}}
\newcommand{\Eb}{\mathbb{E}}
\newcommand{\C}{{{\mathbb C}}}

\newcommand{\U}{{\mathbb U}}

\newcommand{\UU}{{\mathbb U}}

\newcommand{\Sb}{{\mathbb S}}
\newcommand{\llbracket}{[[}
\newcommand{\rrbracket}{]]}
\newcommand{\n}{{\mathcal N}}

\newcommand{\Wc}{\mathcal W}
\newcommand{\one}{\mathbf{1}}

\emergencystretch=\maxdimen
\newtheorem{Theorem}{Theorem}

\newtheorem {lemma} [Theorem]    {Lemma}

\newtheorem {proposition}[Theorem]    {Proposition}

\newtheorem {theorem}[Theorem]    {Theorem}

\newtheorem{remark}[Theorem]{Remark}
\newtheorem{resulta}[]{Result A.}
\newtheorem{resultb}[]{Result B.}

\newcommand {\eps}{\varepsilon}

\usepackage[colorlinks=true,  linkcolor=blue, citecolor=blue,  urlcolor=hair!80!black]{hyperref}

\title[{Parity questions in critical planar Brownian loop-soups}]{Parity questions in critical planar Brownian loop-soups (or ``where did the free planar bosons go?'')}
\author{Matthis Lehmkuehler}
\address {ML: University of Basel --
Departement Mathematik und Informatik}
\author {Wei Qian}
\address {WQ: City University of Hong Kong --
Department of Mathematics and Department of Physics (on leave from CNRS, Universit\'e Paris-Saclay, Laboratoire de Math\'ematiques d'Orsay)}
\author {Wendelin Werner}
\address {WW: University of Cambridge -- Department of Pure Mathematics and Mathematical Statistics}

\begin{document}

\begin {abstract}
The critical two-dimensional Brownian loop-soup is an infinite collection of non-interacting Brownian loops in a planar domain that possesses some combinatorial features related to the notion of ``indistinguishability'' of bosons. The properly renormalized occupation time field of this collection of loops is known to be distributed like the properly defined square of a Gaussian free field. In the present paper, we study how much information these fields provide about the loop-soup. Among other things, we show that the exact set of points that are actually visited by some loops in the loop-soup is not determined by these fields.  We further prove that given the fields, a dense family of special points will each have a conditional probability $1/2$ of being part of the loop-soup. We also exhibit another instance where the possible decompositions (given the field) into individual loops and excursions can be grouped into two clearly different groups, each having a conditional probability $1/2$ of occurring.
\end {abstract}

\maketitle 
\section {Introduction}

In this paper, we explore aspects of the relation between the critical Brownian loop-soup in a planar domain and its (renormalized) occupation time field (which is known to be distributed as the square of a Gaussian free field). Let us first recall some definitions and known facts:

The Brownian loop-soup (as introduced in \cite {MR2045953}) in an open planar domain $D$ (with non-polar boundary) is a random countable collection of Brownian loops $(\lambda_i)_{i \in I}$ in $D$. Each Brownian loop $\lambda_i$ has a time length $T_i$ but no specific marked root on it. When $D$ is bounded, then for each $\eps >0$, the number of loops with time-length greater than $\eps$ is finite; on the other hand, there are infinitely many small loops, and these small loops are dense in $D$. The Brownian loop-soup is a Poissonian collection, which loosely speaking means that the Brownian loops appear ``independently'' from each other. The {\em critical} Brownian loop-soup obtained for one specific intensity of loops is of particular interest, as it possesses a {\em rewiring property} \cite {MR3618142} that can be related to the notion of {\em indistinguishability of bosons} (loosely speaking, this property means that the law of the loop-soup is invariant under a natural Markovian process, where two different overlapping loops can be concatenated into one loop and conversely loops can be cut into two smaller ones at double points).

This rewiring property actually leads naturally to focus on the loop-soup clusters (introduced in \cite {MR2979861}) defined by such a loop-soup:  Two loops $\lambda$ and $\lambda'$ in a Brownian loop-soup $(\lambda_i)_{i \in I}$ are said to be part of the same cluster if it is possible to find a finite chain $i_0, \ldots, i_n$ in $I$, such that $\lambda = \lambda_{i_0}$, $\lambda' = \lambda_{i_n}$, and $\lambda_{i_j} \cap \lambda_{i_{j-1}} \not= \emptyset$ for all $j \in \{ 1, \ldots , n\}$. Intuitively, the aforementioned rewiring property means that the conditional law of the decomposition of the cluster into individual loops is ``uniform'' (and this can be made precise in the discrete settings, see again \cite {MR3618142}).

As we shall recall in more detail in Section \ref {S2.1}, these loop-soup clusters are in fact very closely related to the Gaussian free field  (respectively its square) that can naturally  {be} coupled with (resp. constructed from) a critical loop-soup -- and indeed, the Gaussian free field (we will use the acronym GFF in the sequel) is sometimes called the {\em bosonic free field} in the physics literature. To understand one of the main points of our paper, it is useful to bear in mind the difference between a cluster of Brownian loops and its closure: The latter will contain points that are ends of infinite chains of smaller and smaller overlapping loops, so that some points in this closure do actually belong to no Brownian loop (and therefore not to the cluster either).

We define the {\em trace of the loop-soup} to be the union of all the loops in the loop-soup. The trace encapsulates exactly the same information as the collection of loop-soup clusters, as these clusters are the connected components of the trace.

\medbreak
We are now ready to turn to the actual content of the present paper. The main results will be properly stated as Theorems~\ref{1/2} and \ref{mainthm} in Sections~\ref{S2} and~\ref{S3}, and in this introduction we will describe them and their consequences in more loose terms as ``Results''.  A consequence of Theorem~\ref{1/2} will be that:
\begin {resulta}
One can couple two critical Brownian loop-soups in $D$ in such a way that:
\begin {itemize}
 \item Their traces are almost surely different:
 There will actually be points in the loops of the first loop-soup that will belong to none of the loops of the second one (and vice-versa).
 \item But their respective collections of closures of loop-soup clusters are almost surely the same.
\end {itemize}
\end {resulta}
One can recall that {at the critical intensity} the knowledge of the closures of the loop-soup clusters is sufficient in order to construct the square of the GFF (and the occupation time measure) -- see \cite {aru2023excursion} and the references therein (see also \cite{JLQ_loopsoup} for a different approach). So, in some way, these closures are the really relevant objects when one makes the coupling with the GFF. If one conversely knows the GFF, then one can recover (via the nested CLE$_4$ picture, see e.g. \cite {MR3994105}) the closures of the loop-soup clusters. So, one can reformulate our previously stated result as follows:
\begin {resulta}
The occupation time field of the {critical} loop-soup (or the GFF that it is coupled with, or its square) does not determine the trace of the loop-soup.
\end {resulta}

This may come as a surprise, since such a statement does not hold in the discrete (or cable-graph) setting. Furthermore, the properties of the GFF seemed naturally related to the aforementioned ``rewiring property'', so that it was in fact natural to conjecture that the fields would enable to reconstruct the communication classes of loop-soups for the rewiring Markov chain, and the present result disproves this conjecture. See Section \ref {S2.3} for more details and comments.

In fact, Theorem \ref {1/2} that we attempt to summarise in Figure \ref {fig:firstmain.pdf} is a more concrete statement:
\begin{figure}[h!]
  \centering
  \includegraphics[width=.6\textwidth]{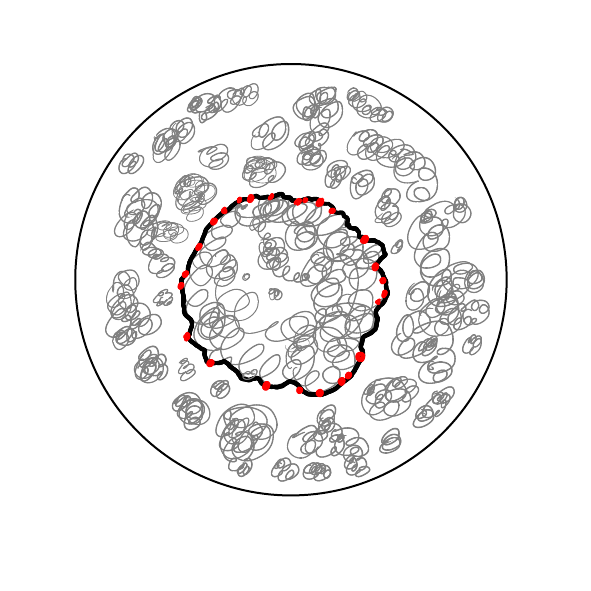}
  \caption{\label{fig:firstmain.pdf}Sketch of our first main result: Conditioning on the excursions by loops away from $\delta$, a whole fractal dense set of points on $\delta$ has a conditional probability $1/2$ of being in the trace of the Brownian loop-soup, and a conditional probability $1/2$ of being in the complement of the trace.}
  \end{figure}
  Consider a critical Brownian loop-soup in the unit disk, and define $\delta$ to be the outer boundary (which is known to be an SLE$_4$-type loop) of the outermost loop-soup cluster that surrounds the origin. It turns out that some points on $\delta$ are part of loops in the loop-soup, but some other points are not part of any loop (i.e. they are in the closure of the cluster but not in the cluster). We then condition on $\delta$ and on the knowledge of all the excursions away from $\delta$ by all the Brownian loops that do touch $\delta$ -- let us call ${\mathcal G}$ this $\sigma$-field.
  \begin {resulta}
We construct two ${\mathcal G}$-measurable dense fractal sets $A_1$ and $A_2$  { of } points on $\delta$, such that conditionally on ${\mathcal G}$:
\begin {itemize}
 \item With conditional probability $1/2$, all points of $A_1$ are in the loop-soup and no point of $A_2$ is in the loop-soup.
 \item With conditional probability $1/2$, all points of $A_2$ are in the loop-soup and no point of $A_1$ is in the loop-soup.
\end {itemize}
\end {resulta}
Let us stress again that (as opposed to many features in the relation between loop-soups and the GFF) this ``quantisation'' phenomenon for the decomposition of closure of Brownian loop-soups into loops is specific to the continuum world, and does not appear in the discrete or cable-system settings.

\medbreak
We now turn to the second part of the paper, that has a similar flavour but will be derived somewhat independently:
It is well-known that it is very natural to superimpose on top of a critical Brownian loop-soup in a domain $D$ an independent Poissonian collection of Brownian boundary-to-boundary excursions inside $D$. Indeed, Dynkin's isomorphism theory (as initiated in \cite {dynkin1983markov,dynkin1984gaussian}) applies to this setup, and (this is not unrelated) it is possible to generalize the aforementioned rewiring property to the union of the collection of excursions and of loops.

Suppose for instance that $R$ is a rectangle of any given width. We consider the union of a critical Brownian loop-soup $\Gamma$ in $R$, and of an independent Poisson point process $\Lambda$ of Brownian excursions in $R$ with endpoints on the vertical sides of $R$ (so either the two endpoints are on the same vertical side, or the excursion crosses the rectangle horizontally -- for such a Poisson point process, there will be infinitely many small excursions with both endpoints on the same side, but only finitely many ones with one endpoint on each side that we call left-right crossings). The number of left-right crossing excursions can change under local rewiring (for instance, when two left-right crossing excursions do intersect, the obtained excursions after rewiring might not be left-right crossings anymore), but the parity of this number of left-right crossings will always be unchanged.
One consequence of our results in Section \ref {S3} will however be that for some excursion intensity,
the parity of the number of left-right crossing excursions is {not} a function of the renormalized occupation field of the union of the excursions and the loop-soup.

\begin{figure}[h]
  \centering
  \includegraphics[width=.6\textwidth]{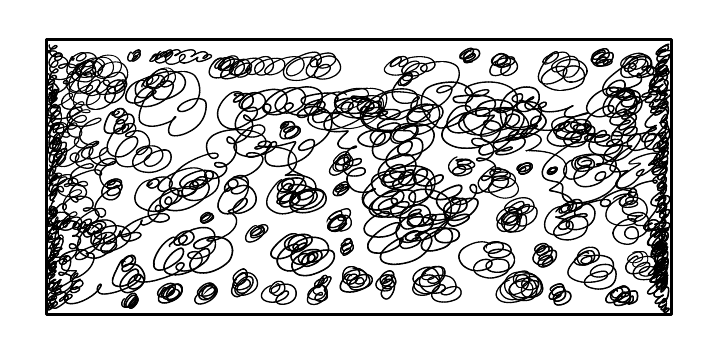}
  \includegraphics[width=.6\textwidth]{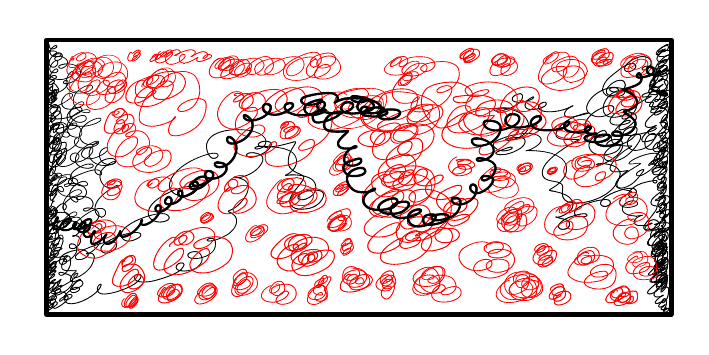}
  \caption{\label{fig:rect}Sketch for Result B.1: The conditional probability that the number of left-right crossing excursions is odd is $1/2$.}
  \end{figure}

  We will in fact show a stronger statement: Let ${\mathcal H}$ now denote the $\sigma$-field generated by the renormalized occupation time field of the union of the excursions and the loop-soup, and let $A$ be the event that a connected component of this union touches both the left and the right-hand side of the rectangle (note that this occurs automatically if there is a left-right crossing excursion, i.e., when $A$ does not occur, then the number of left-right excursions has to be $0$, but $A$ can also occur when the number of left-right excursions is $0$). Then:

\begin {resultb}
On the event $A$, the conditional probability given ${\mathcal H}$ that the number of left-right crossings is odd is almost surely equal to $1/2$.
\end {resultb}

In other words, if one knows the occupation time measure and knows that one cluster touches both sides of the rectangle, then the conditional probability that the number of left-right crossing excursions is odd is always $1/2$.

This result will be closely related to our description of the conditional law of partially explored loop-soup cluster boundaries via parity-constrained Poisson point processes of excursions that we now just heuristically describe -- we refer to Theorem \ref {mainthm} for the precise general statement. Suppose that one starts discovering a {critical} Brownian loop soup in a disk starting from $n \ge 2$ different points (for instance, by following  lines from the boundary as shown in Figure \ref{fig:SecondMain} in the case $n=5$). Each of these $n$ explorations consists of going around (for instance in  counterclockwise manner) the boundaries of all the loop-soup clusters that one encounters (and in their order of appearance on the radial line), and to stop while exploring one of these boundaries. Then, we have $n$ parts $\xi_1, \ldots, \xi_n$ of cluster boundaries which have not been completely explored. Let us consider all the Brownian loops touching $\cup_{i \le n} \xi_i$ and decompose each of them into their excursions away from $\cup_{i \le n} \xi_i$ (note that these will necessarily be excursion from the ``inside sides'' of the $\xi_i$).  {Note that  for each $i$,} the number of excursions starting at $\xi_i$ and ending at $\cup_{j \not= i} \xi_j$ has to be even. We call this {\em  the parity constraint}. Loosely speaking, Theorem \ref{mainthm} then says that:

\begin{figure}[h!]
  \centering
  \includegraphics[width=.6\textwidth]{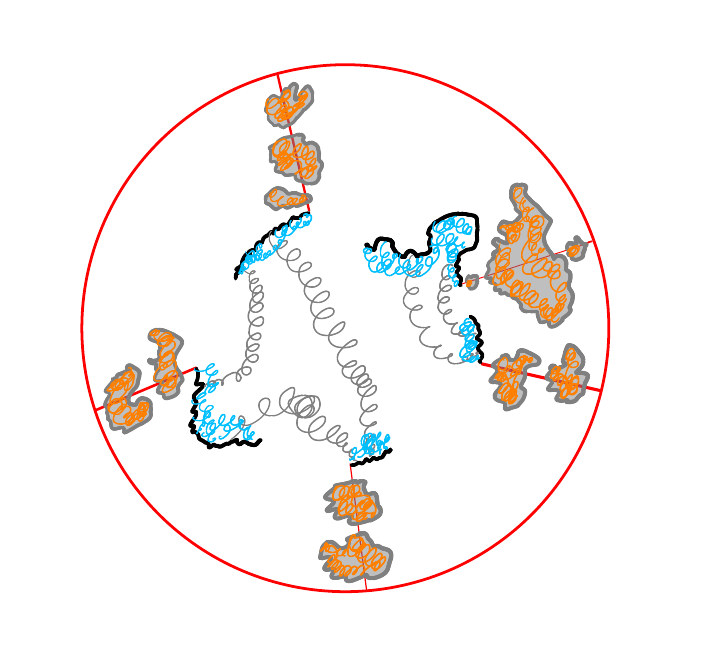}
  \caption{\label{fig:SecondMain}Sketch for Result B.2: The excursions away from $\cup_{1\le i \le 5} \xi_i$ form a Poisson point process {\em conditioned on the parity constraint} (which is satisfied in this example).}
  \end{figure}

\begin {resultb}
 Conditionally on this $n$-strand exploration, the collection of excursions away from $\cup_{i \le n} \xi_i$ is a Poisson point process of Brownian excursions (in the remaining to be explored domain) conditioned to satisfy the parity constraint.
\end {resultb}

The paper is structured as follows: In Section \ref {S2}, we state and prove Theorem \ref {1/2}, and discuss some of its consequences. In Section \ref {S3}, we state and prove Theorem \ref {mainthm} and {also} discuss some of its consequences.
In both cases, the relation between loop-soup clusters, Gaussian free fields and conformal loop ensembles will be key -- we will build on Dynkin's isomorphism type results and earlier papers that we review at the beginning of each sections.

\section{Reconnecting excursions into loops}
\label {S2}

\subsection {Review on loop-soups and the GFF}
\label {S2.1}

We now briefly recall aspects of the relation between the critical Brownian loop-soup and the GFF. As in the introduction, we consider a simply connected domain $D$ with non-polar boundary, and a critical Brownian loop-soup $(\lambda_i)_{i \in I}$ in $D$.
The (renormalized) occupation time measure of such Brownian loop-soups is of particular interest. Let us briefly recall how it is defined:
For each Borel set $B\subseteq D$, one can define $T_i (B)$ to be the time spent by the loop $\lambda_i$ in $B$. As it turns out, for every open set $B$, the quantity $\sum_{i} T_i (B)$ is infinite, because of the many small loops, but one can nevertheless define a renormalized total occupation time ${\mathcal T}(B)$ of $B$ as the limit (in $L^2$ for instance) as $\eps \to 0$ of ${\mathcal T}_\eps (B) - \Eb [ {\mathcal T}_\eps (B) ]$, where ${\mathcal T}_\eps (B) := \sum_{i\in I :\ T_i \ge \eps} T_i (B)$, whenever $B$ is bounded.
In other words, one first forgets about the loops of time-length smaller than $\eps$, recenters the obtained total occupation field, and then takes the $\eps \to 0$ limit. The collection of random variables $({\mathcal T}(B))$ is the {\em renormalized occupation time field} of the loop-soup.
Then (see \cite {le2010markov, MR2815763}), one can see that this field has the same distribution as what is known as the square of a Gaussian Free Field (with Dirichlet boundary conditions) in $D$. Again, if one starts from a GFF, one can construct its square via a  {renormalising/recentering} procedure (so that this obtained ``square'' is in fact also a centered random process).

It is also possible to actually construct the GFF (and not just its square) starting from a critical loop-soup. For this, the notion of {\em loop-soup clusters} that we recalled in the introduction turns out to be instrumental. If $(C_j, j \in J)$ denotes the collection of all these clusters, one can (deterministically) associate a measure $\mu_j$ supported on the closure $\overline C_j$ of $C_j$, and show that $\sum_{j} \Xi_j \mu_j$ is actually a GFF when $\Xi_j$ are chosen to be i.i.d. uniform variables in $\{ -1, +1 \}$, and that the square of this GFF is indeed the occupation time measure of the loop-soup that one started with (see \cite {aru2023excursion,MR3994105,JLQ_loopsoup} and the references therein).

It is probably useful here to make a few historical remarks here:
(1) This Brownian loop-soup can arguably be traced back at least in some implicit way to the early days of Euclidean field theory, where sums over discrete random walk loops were considered by Symanzik (see for instance \cite {Symanzik}), which then led to the construction by several authors of fields as gases of {\em interacting} loops (the goal was indeed to construct non-Gaussian fields) -- this motivated also the definition and study of renormalized self-intersection local times for Brownian motions initiated in \cite {Varadhan}.
(2) This relation between the square GFF and loop-soup occupation times can be considered as a variant of Dynkin's \emph{isomorphism theory} that relates occupation times of Markov processes to Gaussian processes, kicked off in \cite {dynkin1983markov,dynkin1984gaussian} and was subsequently developed in many works (see for instance \cite{eisenbaum2000ray,marcus2006markov} and the references therein).
(3) These results have simple analogue in the discrete setting, when the domain $D$ is replaced by a finite graph. In that case, no renormalisation is needed and one has an identity between non-negative fields (the occupation times on the one hand, and the square of the discrete GFF on the other hand), see \cite {le2010markov,MR4466634} and the references therein.
In fact, many results in the continuum do build on the cable-graph approach to these questions developed by Lupu \cite {MR3502602}.

 {
\subsection {CLE and GFF explorations}
The relation between loop-soup occupation times and the GFF that we described in the previous section lies  at the core of the present paper. We will also build on a number of other inputs, such as the Schramm-Sheffield GFF/SLE$_4$ coupling \cite {MR3101840} and its generalisation to the GFF/CLE$_4$ coupling \cite {MS}, the construction and characterization of Conformal loop ensembles \cite {MR2979861}, Lupu's work on cable graph loop-soups and its consequences \cite {MR3502602} (allowing to match computations on cable graphs with those from \cite {MR3035764}), and of course the paper \cite {MR3994105} in which this type of question about decompositions of loop-soup clusters in the continuum was first addressed (and we will review some of its results relevant for the present paper in the next section).
The relation between loop-soups, CLE$_4$ and the GFF has been used and further studied in several recent papers, including \cite {MR4399157,MR4574830}.}

{
It is maybe worth recalling here one fundamental building block that many of the aforementioned papers build upon and that will be used on several instances in the present paper, namely the ``Markovian explorations'' of Conformal Loop Ensembles that were introduced in \cite {MR2979861}. Suppose that one has a {critical} loop-soup in a domain $D$. The collection $(\delta_i)_{i \in I}$ of its outermost cluster-boundaries form a non-nested CLE$_4$. Choose now any simple curve $L$ starting from the boundary of $D$. Moving along this curve, one will encounter in an ordered way some of these CLE loops. When $L$ bounces into a loop $\delta_i$, one can think of it discovering the whole loop $\delta_i$ ``at once'' (and this is then related to the idea made rigorous in \cite {MR2979861} that this discovery of loops is related to a Poisson point process of SLE-bubbles) or one can choose to discover $\delta_i$ ``progressively'' by going around it continuously in clockwise of anticlockwise manner (and this lies at the core of the proof of the fact that these loops are indeed of SLE-type, as proved in \cite {MR2979861}). In both cases, one important feature (that we will detail in the next section) is then that conditionally on $\delta_i$ (when discovered in this way ``from the outside''), the structure of the cluster of {Brownian} loops that has $\delta_i$ as its outer boundary will depend on $\delta_i$ only in a ``conformally invariant'' way.}

{
The second main feature is that the previous Markovian explorations of the {CLE$_4$} in fact exactly correspond to Markovian explorations of the GFF that it is coupled with, {as implied by a simultaneous coupling of the Brownian loop soup, CLE$_4$ and the GFF established in \cite{MR3994105}.} The loops $\delta_i$ are then the outermost level-lines (with height-jumps $\pm 2 \lambda$) in the CLE$_4$ description of the GFF.  In this setting, the discovery can therefore contain the additional information about the sign of the height-jump along the loops, and the discovered sets are then local sets of the GFF. This will be of particular importance in Section \ref {S3}.}

\subsection {Summary of some previous results and statement of Theorem \ref {1/2}}
\label {Section reconnection} 
\label {S2.2}

Let us first present a survey of some specific previously derived results (in particular from \cite {MR3994105}, we will also try to use notations similar to those of that paper). In this section, we will  consider a {critical Brownian loop-soup} $\Lambda$ in a simply connected planar domain $D  \not= \C$. By conformal invariance, we can for instance assume that $D$ is the unit disc $\UU$. This is the loop-soup with intensity $c=1$ (in the notations of \cite {MR2979861}) or $\alpha=1/2$ (with the notations of \cite {MR3502602}). It has been proved in \cite {MR2979861} that the collection of outer boundaries of outermost loop-soup clusters of this loop-soup form a (non-nested version of a) Conformal Loop Ensemble CLE$_4$ in $D$, that can be alternatively constructed by SLE-type curves as described in \cite {MR2494457} (and a number of features are then accessible via SLE means, such as the fractal dimensions of various sets involved -- see some examples in Section \ref {dimensions}).
Furthermore and more importantly for the present paper, given the collection of these outermost boundaries, the conditional law of the Brownian loops that are surrounded by each of these CLE$_4$ loops turns out to be conformally invariant \cite {MR3994105}. More specifically, one can for instance fix a given point $z$ in $D$, and consider the CLE$_4$ loop $\gamma_z$ that surrounds $z$ (so, $\gamma_z$ will be equal to one of the $\delta_i$'s). This is a simple continuous loop, so that {by Carath\'eodory's Theorem}, the unique conformal map $\psi_z$ from the inside $C_z$ of this loop onto the $\UU$ with $\psi_z (x) = 0$ and $\psi_z' (z) \in \R_+$ does extend continuously to a bijection from $\overline C_z = C_z \cup \gamma_z$ into $\overline \UU$. So, if $z$ was chosen to be the origin, then in Figure \ref {fig:firstmain.pdf}, the loop $\gamma_0$ is depicted in bold, and the map $\psi_0$ sends the region $\overline C_0$ surrounded by $\gamma_0$ (including its boundary $\gamma_0$) into the closed unit disc.  Then, it is shown in   \cite {MR3994105} that:

\begin {enumerate}
 \item
 If we define by $\Lambda_z$ the collection of Brownian loops that are contained in $\overline C_z$, then $\tilde \Lambda_z := \psi_z ( {\Lambda}_z )$ is independent of $\gamma_z$.
\item The law of $\tilde \Lambda_z$ is invariant under any fixed Moebius transformation of $\UU$ onto itself.
\end {enumerate}
By conformal invariance of the Brownian loop-soup itself, the collection $\tilde \Lambda_z$ is in fact also independent of $D$ and $x$. From now on, we will work with this random collection of loops in the unit disc -- that we will call by $\tilde \Lambda$.  {The point $z$ and the original domain $D$ in which the loop-soup was sampled will be fixed to be the origin and the unit disk, and we will just write $\psi$ for $\psi_z$ (we will now also use $z$ to denote other points).}

Further results derived in \cite {MR3994105,MR3896865} (using a combination of ideas, in particular the close relation between the occupation times of this loop-soup and the square of the Gaussian Free Field) include the following:
\begin {enumerate}
\setcounter{enumi}{2}
\item
 If we split $\tilde \Lambda$ into two parts $\tilde \Lambda^b$ and $\tilde \Lambda^i$  defined as the collection of loops of $\tilde \Lambda$ that do touch the boundary $\partial \UU$ and the collection of loops that do not touch $\partial \UU$, {respectively,} then these two sets of loops are independent, and $\tilde \Lambda^i$ is just a critical loop-soup in $\UU$.
\item  { If we split each loop of $\tilde \Lambda^b$ into its countably many excursions away from $\partial \UU$, then one can consider the collection $\tilde \Sigma = (\tilde \sigma_i)_{i \in I}$ of all excursions made by all loops in $\tilde \Lambda^b$. This $\tilde \Sigma$ in now a countable collection of excursions away from $\partial \UU$ in $\UU$. This set then turns out to be a Poisson point process of Brownian excursions with an intensity $\mu / 4$, where $\mu$ is the Brownian excursion measure in $\UU$ (throughout this paper, we use the following normalization for our Brownian excursion measures: We first normalize the Brownian excursion measure in the upper half plane so that the
corresponding measure on pair of endpoints $u,v$ on the real line has a density $(u-v)^{-2} du dv$ -- the excursion measure in any other simply connected domain is the defined from this one via conformal invariance).}  This is in fact proved in two steps: First, it is shown {in \cite{MR3994105}} that the laws of the occupation time measures of the union of the excursions are the same for these two collections of excursions. Then in \cite {MR3896865}, it is proved that the law of its occupation time measure does indeed determine the law of the point process of excursions.
\end {enumerate}
To fully complete the picture, one would have to understand how to construct the law of $\tilde \Lambda^b$ out of the collection of excursions $\tilde \Sigma$, i.e., how to wire these Brownian excursions back into the collection of boundary-touching loops.

Let us now describe the heuristic argument presented in \cite {MR3994105} that suggests that $\tilde \Lambda^b$ is {\em not} a deterministic function of $\tilde \Sigma$. For each $x \in \UU$, define $W_x$ to be the connected component that contains $x$ of $\UU \setminus \cup_{i \in I} \tilde \sigma_i$.
Since the intensity of the excursion measure in the Poisson point process is $\mu /4$, the union of the excursions are related to restriction measures of parameter $1/4$ which are described by SLE$_{8/3} ( \rho)$ processes with $\rho=(-8+2\sqrt{7})/3$. Therefore the event ${\mathcal W}_x := \{ \partial W_x \cap  \partial \UU \not= \emptyset \}$ has positive probability. On $\Wc_x$, let $I_x$ be the set $\partial W_x \cap  \partial \UU $. It is easy to see (for instance, using the Markov property of the Bessel process that drives this SLE-type process) that on the event that $I_x$ is not empty, then it has almost surely no isolated point. Some points in $I_x$ can be isolated from one side (on the unit circle), but there are at most countably many such points. On the other hand, when $I_x$ is not empty, then the Hausdorff dimension of $I_x$ is some positive constant (see Section \ref {dimensions} for a brief discussion of the actual fractal dimensions). One can also note (for instance via a simple $0-1$ law argument looking at the $\sigma$-field generated by small excursions, and combine this with conformal invariance) that $\cup_{x \in \UU} I_x$ is almost surely not empty and dense on $\partial \UU$.

\begin{figure}[h]
  \centering
  \includegraphics[width=.8\textwidth]{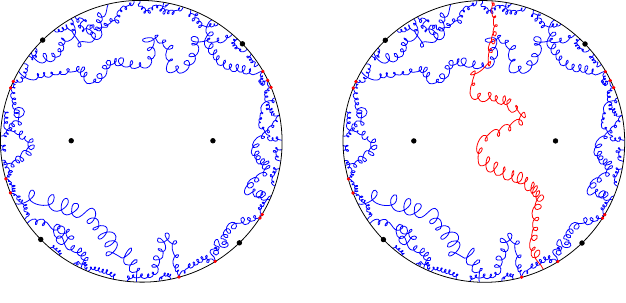}
  \caption{\label{fig:loop-soup-extra}Both the left and the right side of the figure show the boundary-to-boundary excursions of the loops of a Brownian loop soup on the inside of a CLE loop. On the left, we consider the event $\mathcal{W}_{1/2}\cap \{W_{1/2}=W_{-1/2}\}$ and marked points on the boundary show the set $W_{1/2}\cap \partial \UU$, while the right side differs by an extra excursion from the bottom to the top as explained in the main text.}
  \end{figure}

Now, when ${\mathcal W}_x$ happens, then:
\begin {enumerate}
 \item
Either no loop in $\tilde \Lambda^b$ surrounds $x$.  In this case, no point in $I_x$ (except possibly the ones that are isolated on one side) belongs to a loop in $\Lambda$.
\item Or there exists exactly one loop in $\tilde \Lambda^b$ that goes around $x$ (and this loop visits each point in $I_x$ exactly once, and no other loop visits any other point of $I_x$).
\end {enumerate}
This (i.e., that if no loop surrounds $x$ then almost all points in $I_x$ do not belong to a loop, that if one loop surrounds $x$, then it visits each point of $I_x$ exactly once and no other loops visits a point in $I_x$, and that it is not possible that two loops surround $x$) follows from the fact that almost surely, Brownian paths have no points that are simultaneously double points and local cut points \cite {MR1062056} (this follows from the fact that the Brownian intersection exponents $\xi (2,2)$ and $\xi (3,1)$ are strictly larger than $\xi (2,1) = 2$, as shown by earlier pre-SLE work of Lawler \cite {MR1117680} --   {for the definition of intersection/disconnection exponents and their relation to the dimension of subsets of the Brownian path, the reader can for instance consult} \cite {Law,LSW1}).

The following argument from \cite {MR3994105} indicates that each one of the two  scenarios (1) and (2) is in fact possible: Suppose that we have a configuration for  {$\tilde \Sigma$} of the type depicted in Figure \ref{fig:loop-soup-extra}. Then (by resampling the set of macroscopic excursions from the top quarter to the bottom quarter of the circle), it can  happen with positive probability that there is a single additional excursion $e_+$ (instead of none) {that separates $-1/2$ and $1/2$} as indicated on the picture. Then, the loop that contains this excursion will either visit all the points of $I_0$ that are in the right quarter-circle (and in that case, the point ${1/2}$ will be surrounded by a loop while ${\mathcal W}_{1/2}$ occurs, and ${-1/2}$ is surrounded by no loop while ${\mathcal W}_{-1/2}$ happens), or the loop will visit all the points of $I_0$ that are in the left quarter-circle (and in that case, the converse events holds) -- again this is simply due to the non-existence of double cut-points on Brownian paths.

The result that we will now prove is the following stronger fact: 
\begin {theorem}
\label {1/2}
On the event ${\mathcal W}_0 = {\mathcal W}_0 (\tilde \Sigma)$, the conditional probability given $\tilde \Sigma$ that no loop in $\tilde \Lambda^b$ surrounds $0$  is almost surely equal to $1/2$ (and therefore the conditional probability that a loop in $\tilde \Lambda^b$ does surround $0$ and goes through all the points of $\partial W_0 \cap \partial \UU$ is also almost surely equal to $1/2$ on ${\mathcal W}_0$).
\end {theorem}

\subsection {Proof of Theorem \ref {1/2}}

It will be more convenient to work in the horizontal strip $\Sb:= \R \times (-1, 1) \subset \C$ instead of $\UU$. We let $\Theta = ( \theta_j)_{j \in J}$ and $\Sigma = ( \sigma_i)_{i \in I}$ denote respectively the conformal images of $\tilde \Lambda^b$ and $\tilde \Sigma$ under some fixed conformal map from $\UU$ to ${\mathbb S}$. So, $\Sigma$ is a Poisson point process of excursions in ${\mathbb S}$.
One can note that its law is invariant under horizontal translations.  Note that all excursions $\sigma_i$ are bounded. The conformal invariance of the configuration {\em of loops} $\tilde \Lambda$ in $\UU$ actually also directly implies that the law of their conformal image in ${\mathbb S}$ is invariant under horizontal translations.

For each point $x \in \Sb$ (we will actually only consider $x \in \R$), we define the connected component $V_x$ of $\Sb \setminus \cup_j \theta_j = \Sb \setminus \cup_i \sigma_i$ that contains $x$ and the event ${\mathcal V}_x$ that the boundary of $V_x$ intersects $\partial \Sb$. We then decompose ${\mathcal V}_x$ into the two events ${\mathcal V}_x^+$ and ${\mathcal V}_x^-$ corresponding respectively to the  cases where there exists a loop in $\Theta$ that surrounds $x$ or not.

Note that the set $V_x$ (and therefore the event ${\mathcal V}_x$) is  measurable with respect to the collection of excursions $\Sigma$. Reformulating the theorem in the setting of the strip, we see that our goal is to show that
$$ \Pb [{\mathcal V}_0^+ | \Sigma ] = \Pb [ {\mathcal V}_0^- | \Sigma ] = \frac {1}2 \times \one_{{\mathcal V}_0} $$
almost surely. In other words, we want to show that for any $\sigma (\Sigma)$-measurable set $A$ that is contained in ${\mathcal V}_0$,
$\Pb [ A^+ ] = \Pb [A] /2$
where $A^+ := A \cap {\mathcal V}_0^+$ -- indeed this means that for all $\sigma(\Sigma)$ measurable set $B$, $\Eb [ \one_B \one_{{\mathcal V}_0^+} ] = \Pb [ \one_B \one_{{\mathcal V}_0}/2 ]$, which allows to conclude by definition of conditional expectation.  

Clearly (by monotone convergence letting $n_0 \to \infty$, noting that the diameter of $V_0$ is anyway a finite random variable), it suffices to prove that for any given $n_0$, $\Pb[A^+ ] = \Pb [ A ] / 2$  for any $A \in \sigma (\Sigma)$ that is a subset of ${\mathcal V}_0 \cap \{ V_0 \subset (-n_0, n_0) \times (-1, 1) \}$, which is what we now proceed to do. We can also of course assume  that $\Pb [ A] > 0$. From now on, this event $A$ will be fixed.

Let us fix $\eps>0$ with $\eps < \Pb[A] /2$. Let $N (\Sigma,[z, z+n])$ denote the number of excursions of $\Sigma$ that stay in $[z, z+ n] \times [-1, 1]$ and join the top of the strip to the bottom of the strip. This is a Poisson random variable with mean $a(n)$ that tends to infinity as $n \to \infty$. In particular, one can find $n_1 \ge n_0$ such that for all $n \ge n_1$,
$$ \Pb [ N ( \Sigma, [z, z+n])\hbox { is even}] \in (1/2 - \eps, 1/2 + \eps )$$
and the same estimate for the probability that this number is odd.

For each $x$ on the real line, we define $A_x$ (respectively $A_x^+$) to be the event
that $A$ (resp. $A^+$)  holds for the picture of $\Sigma$ (resp. $\Theta$) shifted  horizontally by $x$ (so that we are looking at the properties of $V_{x}$ instead of $V_0$).

The proof will build on the following observation: If $A_x^+$ and $A_{y}^+$ both hold (or if $A_x^-$ and $A_{y}^-$ both hold), then the number of excursions in $\Sigma$ that disconnect $y$ from $ x$ in $\Sb$ is necessarily even, while if $A_x^+$ and $A_{y}^-$ both hold (or if   $A_x^-$ and $A_{y}^+$ both hold), then this number is necessarily odd. On the other hand, when $x$ and $y$ are more than $2n_0 + n_1$ apart, this number of excursions is a Poisson random variable with very large mean, which has a probability close to $1/2$ to be even.

More specifically, let us define $m := 3n_1$, and let us define $R$ to be the collection of all rectangles $(jm + n_1 , jm+ 2n_1 ) \times (-1, 1)$ for $j\ge 1$. Let ${\mathcal G}$  be the $\sigma$-field generated by the set of  all excursions of $\Sigma$ that do not fully stay in any rectangle in $R$. Then, we first observe that all events $A_{jm}$ are measurable with respect to ${\mathcal G}$.

\begin{figure}[h]
  \centering
  \includegraphics[width=\textwidth]{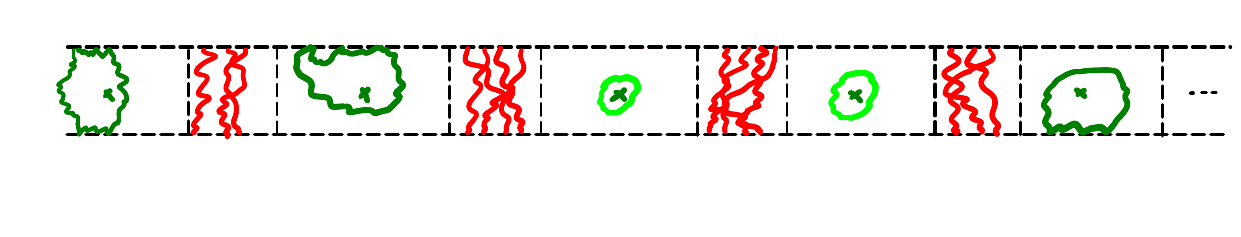}
  \caption{\label{fig:Strip} Idea of the proof: We look at which $A_{jm}$'s hold (here this happens for three out of the five depicted $j$'s), and note that the number of separating top-to-bottom excursions has a probability close to $1/2$ to be even, which implies that the proportion of such $j$'s for which $A_j^+$ holds is close to $1/2$ (in practice, the boxes will be much wider, since $m$ will be chosen very large in order for the probability that the number of top-to-bottom crossings in a box is even to be close to $1/2$).}
  \end{figure}

Let $J$ denote the set of positive integers for which $A_{jm}$ hold, and we let $J^K := J \cap \{1, \ldots, K\}$. 
Let us first note that $\# J^K / K $ does almost surely tend to $\Pb [A]$. This follows directly from Birkhoff's Ergodic Theorem (see e.g., Theorem 6.21 in \cite {MR1163370}). Indeed, the law of the Poisson point process of excursions is clearly invariant under horizontal shift $T$ by the  $m$, and this transformation $T$ is easily shown to be ergodic.
(For each measurable event $B \in \sigma(\Sigma)$, one can find $B_N$ that is measurable with respect set of excursions in $\Sigma$ that are contained in the box $[-N, N] \times [-1,1]$ so that
$\Pb [B \Delta B_N] \to 0$ as $N \to \infty$. Since for each given $N$, the set $T^N [B_N]$ is independent of $B_N$ and if $B$ is invariant under $T$,
$$ \Pb [ B \Delta T^N (B_N) ] = \Pb [ T^N (B) \Delta T^N (B_N)] =
\Pb [ B \Delta B_N ] \to 0, $$
one easily concludes that any $T$-invariant event has probability $0$ or $1$.)

We can in particular always define 
$j_0 := \min J$. 
If we condition on ${\mathcal G}$, then we know $J$ and $j_0$. 
We then call $J_0$ (respectively $J_1$) the subset of $J$ consisting of those $j$ for which $jm$ and $j_0m$ are separated in $\Sb$ by an even (resp. odd) number of excursions of $\Sigma$. We also call $J_0^K$ and $J_1^K$ the respective intersections of $J_0$ and $J_1$ with $\{1, \ldots, K\}$. 
The previous parity observation shows that the number $N_K := \sum_{j =1}^K \one_{A_{jm}^+}$ is equal to one of the two values $\# J_0^K$ or $\# J_1^K$, so that almost surely,  
$$ \min (\# J_0^K, \#J_1^K) \le N_K \le \max ( \#  J_0^K, \# J_1^K ). $$ 
The idea will be to estimate $\Eb [ N_K]$ (which is also equal to $K\Pb [A^+]$ by translation invariance) exploiting the ergodicity of horizontal shifts for the Poisson Point process $\Sigma$ and the fact that $\# J_0^K$ and $\# J_1^K$ are actually functions of $\Sigma$ only.

Let us combine the following observations: 
\begin {itemize} 
 \item 
If we condition on ${\mathcal G}$, we know $J$ and have some information about excursions separating the $jm$'s, but we are still missing the information about the top-to-bottom excursions that lie in $R$. In particular, we see that conditionally on any event in ${\mathcal G}$, we can always find a coupling to upper bound the realizations of the conditional laws of $\# J^K_0 -1$ and $\# J^K_1$ by the sums of $\# J - 1$ independent Bernoulli random variables with mean $1/2 + \eps$.
\item When $u$ is large enough (say $u \ge u_0$), we know that the probability that the sum of $u$ independent Bernoulli random variables with mean $1/2 + \eps$ is  larger than $u (1/2 + 2 \eps) -1$ is smaller than $\eps$. 
\item We know that
$\# J^K  / K \to \Pb [ A ]$ almost surely. In particular, when $K$ is chosen large enough,
$$ \Pb [ \# J^K / K \le \Pb[A] - \eps ] \le \eps .$$
\end {itemize} 
Hence, if we choose $K$ large enough so that $K (\Pb[A]-\eps) > u_0$, we get that 
$$ \Pb [ \# J_0^K \ge K (1/2 + 2 \eps ) \Pb [A] ] \le 2 \eps.$$ 
Combining this with same argument applied to $\# J_1^K$ then yields that for large enough $K$, 
$$ \Pb [ N_K \ge K (1/2 + 2 \eps ) \Pb [A] ] \le 4 \eps.$$ 
The very same type of argument applied to $A^-$ instead of $A^+$ shows on the other hand that for large enough $K$, 
$$ \Pb [ N_K \le K (1/2 - 2 \eps ) \Pb [A] ] \le 4 \eps.$$ 
We now conclude that $N_K / K$ converges in probability to $\Pb[A]/2$ as $K \to \infty$, so that by dominated convergence, $\Eb [ N_K / K ] \to \Pb [A]/2$. But by translation invariance, $\Eb [ N_K / K ] = \Pb [ A^+]$ so that we can conclude.

\subsection {Further comments and results}
\label {S2.3}

\subsubsection {About the rewiring property}
As already mentioned, the critical Brownian loop-soup has a striking and very natural ``rewiring property'' (see \cite{MR3618142}) that in turn can be shown to imply that the law of the loop-soup is in fact invariant under Markov chains that can be loosely speaking described as follows. For each fixed $\eps$, we can define a Markov Chain ${\mathcal M}_\eps$ on the state-space of the Brownian loop-soups as follows (see \cite {MR1229519} for background on intersection local times):

\begin{itemize}
  \item Any two disjoint loops of time-length at least $\eps$ in the loop-soup that do intersect, will merge into a single loop at a rate and place provided by the intersection local time between these two loops.
  \item Any single loop will split into two loops of time-length at least $\eps$, at a rate and position provided by the self-intersection local time on the loop restricted to the set of times $0< s< t \le T_i$ where $t-s > \eps$ and $s+ T_i - t > \eps$ (this just means that the two obtained loops both have time-length at least $\eps$).
\end{itemize}

These Markov chains do leave the occupation time field of the Brownian loop-soup as well as the loop-soup clusters unchanged. This made it tempting to conjecture that this occupation time field (or equivalently, the closure of the loop-soup clusters) does in fact determine the communication class of the loop-soup under these Markovian dynamics (i.e., the union over $\eps$ of these communication classes), which is nothing else than the loop-soup cluster itself.
Indeed, such a statement could be shown to hold in the setup of cable-graphs (with loop-soups on cable-graphs, as studied in \cite {MR3502602}). Theorem \ref {1/2}  shows that this is not the case.

Note that the natural dynamic that would be involved in switching between the two options exhibited by Theorem \ref {1/2} (since both have conditional probability $1/2$, it is easy to update/resample the outcome for each cluster boundary independently without changing the law of the loop-soup) would be very non-local and involve at least ``switching'' the status of uncountably many points from ``visited by a loop'' to ``not-visited by any loop'' and vice-versa.

\subsubsection {The results A.1, A.2, A.3 stated in the introduction}

A first comment is that it is possible to use Theorem \ref {1/2} for other related cases than just the outer boundaries of the outermost loop-soup cluster containing the origin. One can for instance:
\begin {itemize}
 \item Use the inversion invariance and consider the inner boundaries of loop-soup clusters instead of the outer boundaries. 
 \item One can first apply the restriction property of the loop-soup, looking at the set of loops that are contained in a given subset $D'$ of the domain $D$, and then look at the special points on the outer boundaries of loop-soup clusters for this smaller loop-soup. One can also use Markovian ways to define $D'$ at random (such as choosing $D'$ to be the inside of a loop-soup cluster).
\end {itemize}
So, there are plenty of possible loop-soup measure preserving switchings that one can do in the same spirit as Theorem \ref {1/2}.

This can be used to check that the set of points which have a conditional probability $1/2$ of being in the trace of the loop-soup is in fact dense in the domain. Indeed, one can note that the set of (non-outermost) loop-soup clusters is dense (as any given point is almost  {surely} surrounded by infinitely many disjoint nested clusters), and there will be such points on the outer boundary of each of these clusters.

\medbreak

Let us now outline how to deduce Results A.1, A.2 and A.3 described in the introduction from Theorem \ref {1/2}. We call ${\mathcal G}$ the $\sigma$-field generated by the information provided by the entire loop-soup in the unit disc, except the knowledge of how the excursions away from $\delta$ (which is the outer boundary of the outermost loop-soup cluster surrounding the origin) are wired together. {In other words, this is the $\sigma$-field generated by the outermost loop-soup cluster boundary $\delta$ surrounding the origin, the Brownian loops that do not intersect $\delta$ and the collection of excursions $\tilde \Sigma$. So, in terms of the image $\tilde \Lambda$ under $\psi$ of the collection $\Lambda$ of Brownian loops surrounded by $\delta$ via $\psi$, one has access to the collections of excursions $\tilde \Sigma$ and to the Brownian loops $\tilde \Lambda^i$ that do not intersect $\partial \U$.}

Let us come back to the notations of Section \ref {S2.2}. 
{For each $x \in \U$, let $I_x$ be the set of points in $\partial W_x\cap \partial \U$ which are not isolated on the left or right (in this set). We remove those points which are isolated on one side, because they could possibly be visited by a loop, even when all points in $I_x$ (which are by definition not isolated on either side) are not visited by any loop. Since there are only countably many points in $\partial W_x\cap \partial \U$ which are isolated on one side, removing them does not change the Hausdorff dimension.
Let $I = \cup_{x \in \U} I_x$. }
Note that any $W_x$ contains points with rational coordinates, so that this union can be taken over the countable set of rational points in $\U$.

Let us now explain why $I$ is almost surely dense on $\partial \U$.
The relation between Poisson point processes of excursions, restriction measures and  SLE$_{8/3} (\rho)$ processes \cite {MR1992830} (see also \cite {MR2178043})  ensures that with positive probability, this  set $I$ is non-empty with positive Hausdorff dimension (because the intensity $1/4$ is strictly smaller than $1/3$ which is the critical one for this question).
We can note that almost surely, $I$ does contain no endpoint of a Brownian excursion, because there are almost surely only countably many such endpoints, and almost surely, each of these endpoints will be overarched by infinitely many small excursions that will prevent it from being in $I$.
The definition of $I$ furthermore immediately implies that for any closed arc $a \subset \partial \U$ with positive length, the set $a \cap I$ is in fact independent of the set of excursions with no endpoint in $a$ (since these excursions will be at positive distance from $a$). Combining these two observations implies that $a \cap I$ is in fact measurable with respect to the set of excursions that have both endpoints in $a$.
By conformal invariance, the probability that $a \cap I$ is not empty does in fact not depend on the length of $a$. So, if one subdivides any subarc of  $\partial \U$ into $N$ disjoint arcs of smaller length (and then letting $N$ to infinity), we can conclude that $I$ is indeed dense on $\partial \U$. We can finally define $A \subset \delta$ to be the image of $I$ under the conformal map $\psi$ from the unit disc back into the interior of the CLE loop.

Let us now turn to the definition of the two sets $A_1$ and $A_2$.
Any given pair of points $y_1$ and $y_2$ on $\partial \U$ does split $\partial \U$ into {two} open boundary arcs. We denote by $N_{y_1,y_2}$ the number of Brownian excursions in the Poisson point process that have one endpoint on each of the boundary arcs (or equivalently that disconnect $y_1$ from $y_2$ in $\U$). On the one hand, for fixed $y_1$ and $y_2$, $N_{y_1, y_2}$ is clearly almost surely infinite (due to the infinitely many small excursions overarching $y_1$), but on the other hand, if $y_1$ and $y_2$ are in $I$, $N_{y_1, y_2}$ is necessarily finite (as otherwise, it would imply that infinitely many small excursions overarch either $y_1$ or $y_2$ which would contradict the fact that these points are in $I$. Furthermore, if $N_{y_1, y_2}$ is even, then $y_1$ and $y_2$ are necessarily of the same type (i.e. either both are on some Brownian loops in the loop-soup, or neither of them is on Brownian loops), and  if $N_{y_1, y_2}$ is odd, then necessarily they are of different type.

We can therefore split $I$ into two disjoint of points $I^1$ and $I^2$ (using some deterministic rule) such that $I^1 \cup I^2 = I$ and $N(y_1, y_2)$ is odd for all $y_1 \in I_1$ and $y_2 \in I_2$. The previous considerations then imply that almost surely, all points in $I^1$ are of the same type, all points in $I^2$ are of the same type, and all points in $I^1$ are of a different type than the points in $I^2$. One can reformulate Theorem \ref {1/2} by saying that each given $y \in I^1$ has a conditional probability $1/2$ to be in some Brownian loop. It therefore follows that with conditional probability $1/2$, all points in $I^1$ are in the trace of the loop-soup and none of the points in $I^2$ is in the trace of the loop-soup, and with conditional probability $1/2$, all points in $I^2$ are in the trace of the loop-soup and none of the points in $I^1$ is in the trace of the loop-soup.
Finally, to check that $I^1$ and $I^2$ are both dense on $\partial \U$, we can note that the proof of Theorem \ref {1/2} (mapping the strip back onto the disc, so that the horizontal half-line $[0, \infty)$ gets mapped to a ray $[0, \exp (i \theta))$) in fact shows that for every given $\theta$, there almost surely exists infinitely many points in $I^1$ and infinitely many points in $I^2$ in the neighborhood of $\exp (i \theta)$. The sets $A_1$ and $A_2$ are then the images of $I^1$ and $I^2$ under $\psi$.

\subsubsection {About fractal dimensions}
\label {dimensions}
The value of the Hausdorff dimensions of the set of points $I$ and $A$ can be derived fairly directly from the
construction and known results. For instance,  the results of \cite {MR1992830} relating restriction measures to SLE$_{8/3} (\rho)$ processes
and the value of the dimension of the intersection of these processes with the boundary as derived in
\cite {MR3602842} (confirming the predictions of \cite {Duplantier}) readily show that the Hausdorff dimension $d_I$ of $I$ is $(\sqrt{7}-1)/6 \sim .27$. The dimension of $A$, which is the preimage of $I$ under $\psi$ can then be obtained via the multifractal spectrum of SLE$_4$ derived in \cite {MR3786302} (it is the maximum of the function $s \mapsto (1+s-2s^2)/ (1-s^2) - (1-d_I)/((1-s))$, which numerically is close to $.285$).

For comparison, let us recall a few known other fractal dimensions here:
The dimension of the set of points on a cluster boundary that do belong to a Brownian loop (and are therefore on the outer boundary of the loop) is in fact equal to $1$ (see \cite {GLQ_loopsoup}, {following from the value of a generalized disconnection exponent \cite{MR4221655}}), which is of course strictly bigger than the dimension of $A$. {It is also proved in \cite{MR4221655, GLQ_loopsoup} that the set of double points (i.e.\ visited at least twice by one loop or at least once by two different loops) on a cluster boundary is $0$.} Recall also that the dimension of the outer boundary of one Brownian loop is $4/3$ (see \cite {MR1849257,MR2153402}), that the dimension of the cluster boundary is $3/2$ (it is an SLE$_{4}$ loop, so that one can use \cite {MR2435854}), and that the  dimension of the set of points surrounded by no CLE$_4$ loop is $15/8$ (see \cite {MR2802511,MR2491617}).

\subsubsection {A natural extension of the set $A$}

Let us now explain how to modify our construction to obtain a set $\hat A$ larger than $A$ with similar properties but larger fractal dimension $(5 - 2 \sqrt {2}) / 4 \sim .54$. This is more a side-comment, so we choose to only outline one natural way to proceed.

The basic general idea is to formally replace each Brownian excursion of $\tilde \Sigma$ in the unit disc by the circular arc (that intersects $\partial \U$ perpendicularly) with the same endpoints. 
We then hook the circular arcs into a collection of loops in the same way as the excursions in $\tilde \Sigma$ are hooked into Brownian loops. The intersection of the loops with the boundary $\partial \U$ remain unchanged, but two arcs will now intersect if and only if their endpoints are intertwined (when the endpoints of the Brownian excursions are intertwined, then they have to intersect, but it can happen that they intersect even if the endpoints are not intertwined), so that the clusters formed by these arcs can be smaller than those formed by the Brownian excursions.

In practice, it is more convenient to work in the upper half-plane instead. So we will work with the map $\varphi$ from the interior of the CLE loop containing the origin into the upper half-plane (instead of the mapping from $\psi$ into the unit disc) that maps the origin on $i$ and has a positive real derivative at the origin, say. One therefore ends up with a Poisson point process of circular arcs in the upper half-plane with endpoints $(u,v)$ chosen with intensity
$\beta du dv / (u-v)^2$ with $\beta=1/4$ (the value $\beta=1/4$ was determined at the end of Section 4 in \cite{MR3994105} -- it comes from the SLE restriction property computation for SLE$_4$ that provides a restriction exponent $1/4$, combined with the fact that a Poisson point process of excursions away from the positive half-line with endpoints chosen with this intensity constructs a one-sided restriction measure of exponent $\beta$ \cite{MR2178043}).
When $y$ is in the upper half-plane, one can look at the connected component $\tilde V_y$ (which is the analog of $W_x$) of the complement of all arcs, and let $\tilde J_y$ be the intersection of $\partial \tilde V_y$ with the real line. Then, for each given $y$, this set $\tilde J_y$ is either empty or a perfect fractal set of positive dimension.

\begin {itemize}
\item We first deduce that the Hausdorff dimension of $\tilde J_y$ is $1/2$ (when it is not empty). This is an easy task and can be done without using Brownian intersection exponents or SLE$_{8/3}(\rho)$ considerations. One can for instance compare it with the set obtained by removing from $\R_+$ a Poisson point process of intervals $(u,v)$ with intensity $\beta  du dv / |u-v|^2$. The obtained set is then the zero-set of a Bessel process with dimension {$4\beta=1$} (which is a standard Brownian motion), i.e.\ a set with dimension $1- 2\beta=1/2$.
The dimension of the preimage of $\tilde J_y$ under $\varphi$ in the loop-soup picture can then be again obtained via the multifractal spectrum of SLE$_4$ derived in \cite {MR3786302} and it turns out to be $(5 - 2 \sqrt {2}) / 4 \sim .54$ (i.e., the maximum of the function $s \mapsto (1+s-2s^2)/ (1-s^2) - 1/(2(1-s))$).

\item
Let $\hat J_y$ be the set obtained by removing from $\tilde J_y$ the countably many points that are isolated on the right, or isolated on the left. Then $\hat J_y$ has the same dimension as $\tilde J_y$.
Suppose that a point $u_0$ in $\hat J_y$ is actually visited by the image under $\varphi$ of one of the Brownian loops. Then, this loop will necessarily have to visit infinitely many other points in $\hat J_y$, and  one of the following two options has to occur (since the loop has to come back to this point $u_0$): (a) That loop will in fact visit all the points of $\hat J_y$ or (b) It visits at least twice all the points of $\hat J_y$ in a right-neighborhood of $u_0$ or in a left-neighborhood of $u_0$.

\item
To conclude, we argue that Scenario (b) is impossible, due to the aforementioned fact that the dimension of double points on the boundary of a loop-soup cluster is $0$ \cite{MR4221655, GLQ_loopsoup}.
In fact, we can also rule out (b) using an earlier result that the dimension of double points on the Brownian frontier (hence the outer boundary of one Brownian loop) is $2- \eta_4 \sim .452 < 1/2$ \cite{MR2644878}, where $\eta_4$ is the disconnection exponent for four Brownian motions \cite{MR1961197}, combined with the fact that the dimension of $\hat J_y$ is greater than $1/2$.
So, either no point in $\hat J_y$ is visited by a loop, or there exists a loop that visits all the points of $\hat J_y$ at least once.
\end {itemize}

From there, we can then proceed as for $I$ and $A$: One can divide $\hat J := \cup_y \hat J_y$ into two sets $\hat J^1$ and $\hat J^2$  with the feature that the conditional probability that every point in $\hat J^1$ is visited once and no point of $\hat J^2$ is visited is $1/2$, and the conditional probability
that every point in $\hat J^2$ is visited once and no point of $\hat J^1$ is visited is $1/2$ as well. The sets $\hat J^1$ and $\hat J^2$ are then dense on the real line  their dimension is almost surely $1/2$.
The dimension of the preimages  $\hat A_1$ and $\hat A_2$ of these sets in the original loop-soup picture is then $(5 - 2 \sqrt {2}) / 4 \sim .54$, which is significantly larger than the dimension of the sets $A_1$ and $A_2$.

 \section {Multiple exploration and the parity condition}
 \label {S3}

 \subsection {Exploring from different points, statement of the main result} 
 We now turn to the results B.1 and B.2 stated in the introduction.
We again consider a critical Brownian loop-soup in the unit disc, which forms a collection of clusters whose outermost outer boundaries form a non-nested CLE$_4$, as shown in \cite {MR2979861}.
Before discussing explorations of the CLE$_4$, let us first recall again from \cite {MR3994105} that if we condition on {\em all} these outermost boundaries $(\delta_j, j \in J)$  (i.e., on the CLE$_4$), then one can divide the set of loops in the loop-soup into the following (conditionally independent) pieces: (a)
  For each $j$, the collection of loops that are surrounded by $\delta_j$ and do not intersect $\delta_j$ form a Brownian loop-soup in the interior of $\delta_j$.
(b) For each $j$, if we look at the collection of all excursions away from $\delta_j$ by those loops that do intersect $\delta_j$, one has a Poisson point process of Brownian excursions in the interior of $\delta_j$, with intensity given by $\mu_j / 4$ (where $\mu_j$ is the standard excursion measure in the interior of $\delta_j$).

Our goal will be to obtain a similar result where one discovers the CLE$_4$ only partially, using
 a Markovian exploration of this CLE$_4$ from $n$ boundary points. Markovian explorations of Conformal Loop Ensembles have been discussed in a number of papers, starting with \cite {MR2979861} itself, and the explorations of CLE$_4$ can also naturally related to the local set theory developed in \cite {MR3101840,MR3477777}.

\begin{figure}[h]
  \centering
  \includegraphics[width=.6\textwidth]{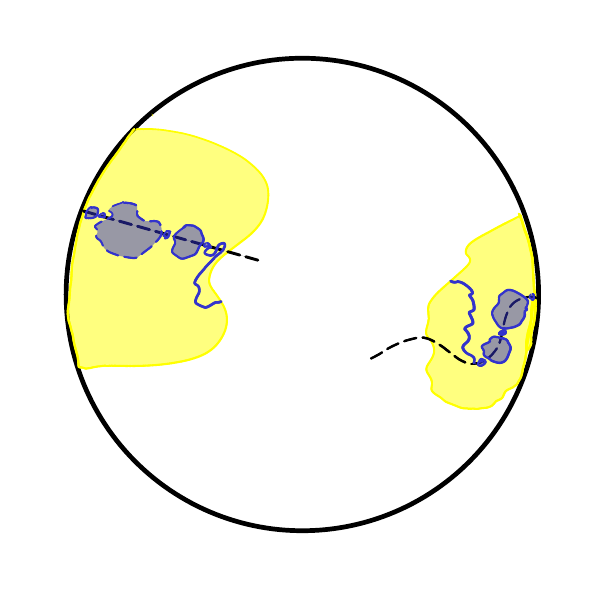}
  \caption{\label{fig:explos} Two stopped explorations obtained by exploring the loops that touch curves $L_1$ and $L_2$ and stopped at the first exit of neighborhoods $O_1$ and $O_2$ of $x_1$ and $x_2$}
  \end{figure}

We are first going to work with some specific explorations (essentially the ones introduced in \cite {MR2979861} in order to characterize the CLEs via their spatial Markov property), and we will then explain later (in Section \ref{MoreGeneral})  how to extend the results to more general explorations.  We choose $n$ given points $x_1, \ldots, x_n$ on $\partial \UU$. From each $x_k$ we grow a straight segment $L_k$ (or another simple curve) towards the inside of the disc, and explore along $L_k$ as depicted in Figures \ref {fig:explos} and ~\ref{fig:SecondMain} (for the $n=2$ or $n=5$ case).
For each $k$, the exploration from $x_k$ traces each loop that $L_k$ encounters in the counterclockwise direction, in the order that $L_k$ encounters them. We then proceed up to some stopping time (with respect to the filtration generated by this exploration along $L_k$) that is defined in such a way that:  (a)  the $k$ explorations remain almost surely  disjoint, and (b) at these stopping times, each exploration is actually in the middle of tracing one of the CLE$_4$ loops -- we call $\xi_k$ the partial piece of this loop that has been traced. One can for instance choose (for each $k$) neighborhoods $O_k$ of $x_k$ in $\UU$ in such a way that $L_k \not\subset O_k$ and the distance between $O_k$ and $O_{k'}$ is positive when $k \not= k'$, and choose to stop the exploration started at $x_k$ at the first time at which it exits $O_k$.

The complement of the exploration starting from $x_k$ is then an open subset $U_k$ of $\UU$ consisting of the inside of the fully-discovered loops and the currently-explored component $U_k'$ which has $\xi_k$ as part of its boundary. We then define $U= \cap_{k=1}^n U_k'$, which is a simply connected subset of the unit disc (recall that the stopped explorations are disjoint by construction). We then let $\varphi$ be a conformal map from $U$ to $\U$ (chosen according to some deterministic rule, among the three-parameter family of such maps).
The ``left-hand side'' of   $\xi_k$ (recall that we trace the loops in the counterclockwise direction, so that this is the ``inside part'' of the CLE$_4$ loop) is mapped to an arc $\partial_k$ on $\partial \U$. The image under $\varphi$ of the rest of the partially explored CLE loops is then distributed as a collection of $n$ disjoint curves (connecting the $2n$ endpoints of the arcs $\partial_1, \ldots, \partial_n$).
Note that there are a Catalan number $C_n$ of possible ways to connect the endpoints of $\partial_1, \ldots, \partial_n$ -- but we will not be discussing these connection probabilities here (see \cite {MW} for results in this direction).

Let us briefly recall two facts from previous papers about the conditional law of the loop-soup given such explorations:
A first fact is the independence of the exploration with respect to the loops that it has not discovered yet -- more precisely,
 the image under $\varphi$ of the loop soup restricted to $U$ is equal to the union of two {\em independent} sets of loops:
\begin{itemize}
\item The set of loops that intersect $\partial:=\partial_1\cup \cdots \cup \partial_n$. This set can be decomposed as a collection $\Sigma_\partial$ of excursions in $\U$ with endpoints in $\partial$.
\item The set of loops that are contained in $\U$. This is distributed as an independent loop soup in $\U$.
\end{itemize}
 {The proof of this fact is similar to \cite[Lemma 4]{MR3994105}:
One has checks that when one resamples the set of loops that are strictly contained in $U$, then the $k$ explorations remain unchanged. As explained in \cite{MR3994105}, there is a slightly subtle feature here: For each given $i$, the Brownian loops that do touch a given $\xi_i$ do not form a unique cluster of loops, but if one adds the Brownian loops of diameter smaller than $\eps$ in the loop-soup to them, then each $\xi_i$ will become automatically part of the boundary of a cluster. From this, one sees readily that for any small $\eps$, resampling the loops of diameter greater than $\eps$ that the explorations haven't hit yet would not change the exploration. Since this is valid for any $\eps$, we conclude that the exploration is indeed independent of the set of loops that it has not hit yet. Note that the result in \cite {MR3994105} is valid for loop-soup of any intensity $c \in (0,1]$, but we are using it here only for $c=1$.}

{
Then, in \cite {MR3901648}, the case of  one ``partial exploration'' (i.e., when $n=1$) was considered for a CLE$_\kappa$ coupled with a loop soup with intensity $c\in(0,1]$. For the critical $c=1$ case, it is {pointed out} there that $\Sigma_\partial$ is a Poisson point process of Brownian excursions away from $\partial = \partial_1$ {(the proof being a rather direct adaptation of the one in \cite{MR3994105})}. More precisely, if $\mu$ denotes the Brownian excursion measure in $\U$ and $E_\partial$ the set of excursions with both endpoints in $\delta$, then the intensity of this Poisson point process is $\mu 1_{E_\partial} / 4$. }

The main goal is now to describe the law of $\Sigma_\partial$ in the cases where $n\ge 2$.
For each $i \not= j$, let $\n_{i,j}$ denote the number of excursions with one endpoint in $\partial_i$ and the other endpoint in $\partial_j$. For each $i$, we then define $\n_i := \sum_{j \not= i} \n_{i,j}$ to be the number of excursions with just one endpoint in $\partial_i$. Note that because of the fact that these excursions are the excursions made by loops,  $\n_1, \ldots, \n_n$ are necessarily all even. So, $\Sigma_\partial$ is clearly not a Poisson point process of excursions away from $\partial$.

A substantial part of the rest of the paper will be devoted to the proof of the following result:

\begin {theorem}
\label {mainthm}
The conditional law of $\Sigma_\partial$ is that of a Poisson point process of Brownian excursions in $\U$ with both end-points in $\partial$ (i.e., with intensity $\mu 1_{E_\partial}/4$)  conditioned by the event that $\n_i$ is even for all $i \in \{ 1, \ldots, n\}$.
\end {theorem}

Note that the event that $\n_i$ is even for all $i \le n$ is in fact independent of the collection of excursions from $\partial$ to $\partial$ such that both endpoints lie in the same $\partial_j$. So, one can restate this theorem by saying that one can decompose $\Sigma_\partial$ into two independent parts: The Poisson point process of Brownian excursions from $\partial$ to $\partial$ with both endpoints in the same $\partial_i$ (and there will be infinitely many such excursions for each $i$), and the point process of excursions with endpoints in different $\partial_i$'s which is a Poisson point process of such excursions conditioned on the fact that $\n_i$ is even for all $i$.

\medbreak
The strategy of the proof will be to exploit the coupling of the CLE$_4$ with the GFF and can be outlined as follows:
\begin {enumerate}
 \item We first couple the loop-soup with a GFF  as described in Proposition 5 of \cite {MR3994105},
  i.e., in such a way that (a) the outermost loop-soup cluster boundaries are the outermost CLE$_4$ level lines defined by the GFF, and (b) the square of that GFF is the renormalised occupation time of the loop-soup. (Throughout this paper, we will choose the normalization of the GFF so that it has generator $-\Delta/2$; this corresponds to the height gap $2\lambda$ in the level line coupling being $\sqrt{\pi}$.)  We then do the previously described partial exploration.
 For each $k$, we define $\sigma_k \in \{-1, +1\}$ to be the value such that $\xi_k$ is part of a $0/2\sigma_i \lambda$ level line of the GFF.
 The first step (Section~\ref{subsec:coupling}) of the proof will be, using GFF partition function type considerations, to determine the joint law of $(\sigma_1, \ldots, \sigma_n)$ conditionally on the exploration as a function of $(\partial_1, \ldots, \partial_n)$.
 \item If one further conditions on the values of $(\sigma_1, \ldots, \sigma_n)$, we get the conditional law of the GFF in $U$, and (via Dynkin's isomorphism) one obtains the conditional law of its square in $U$. By averaging over the possible values of $(\sigma_1, \ldots, \sigma_n)$ using the law of this $n$-tuple derived in the first step, we get the conditional law of the occupation time field given the $n$ partial explorations (Section \ref {S3.3}).
\item We know that conditionally on $U$, the collection of loops that fully stay in $U$ (without touching its boundary) is a standard Brownian loop soup in $U$. As a consequence, we can deduce from the previous step the conditional law of the remaining part, i.e., of the sum of the occupation times of the boundary-touching loops (i.e., after mapping onto $\U$,  corresponding to loops that do touch $\partial_1 \cup \ldots \cup \partial_n$) -- this will still be derived in Section \ref {S3.3}.
\item The next step is to show that this law is the same as that of the occupation time of the parity-constrained Poisson point process of excursions described in the theorem. This will use again Dynkin Theorem type considerations and a splitting rearrangement trick reminiscent of the ones appearing in the random current representation of the Ising model (Sections \ref {S3.4} and \ref {S3.5}).
\item It then finally remains to argue (see Section \ref {S3.6}) that this identity between the laws of the occupation times is sufficient to deduce that the law of the point processes of excursions themselves are indeed the same. This last step will be based on  the fact that the occupation time of a point process of Brownian excursions determines the law of this point process (which has been derived in \cite {MR3896865}).
\end {enumerate}

\subsection {The exploration and the coupling with the GFF} \label{subsec:coupling}

We are going to be in the setup described in the previous section (with $n$ disjoint stopped explorations of the same CLE$_4$, and coupling with the GFF and the labels 
$\sigma_1, \ldots, \sigma_n$). Let us stress that the $\sigma_k$'s are not independent of the remaining-to-be discovered configuration of loops (for instance, if $\sigma_k \not= \sigma_{k'}$, then $\xi_k$ and $\xi_{k'}$ are necessarily not part of the same CLE loop).

For each $i,j \le n$, we define $\mu_{i,j}$ to be the standard Brownian excursion measure in $\U$, restricted to those Brownian excursions that have one endpoint on $\partial_i$ and the other endpoint in $\partial_j$.
The total mass $|\mu_{i,j}|$ of $\mu_{i,j}$ is finite when $i \not= j$, and we define {$m_{i,j} := | \mu_{i,j} | / 4$.}
The goal of this section is to show that:

\begin {proposition} 
\label {mainprop2}
The conditional probability $p(a_1, \ldots, a_n)$ that $(\sigma_1, \ldots, \sigma_n) = (a_1, \ldots, a_n)$ given the $n$ explorations, is (for all $a_1, \ldots, a_n \in \{-1, 1\}$)
\begin{align}\label{eq:prop2}
p( a_1, \ldots, a_n)
= \frac {1}{Z_{{n}}} \exp \Big(  \sum_{i < j}  a_i a_j m_{i,j}\Big)
 \end{align}
where $Z_{{n}}$ is the normalizing factor so that the sum of the probabilities is $1$.
\end {proposition}

Let us give a first heuristic justification/explanation for this result. A first remark is that  if for all $i$, $\Phi_i$ denotes the harmonic function in $\U$ with boundary condition { $ 2 \lambda \one_{\partial_i}$,} then when $i \not= j$,  $2 m_{i,j}$ is equal to
$$(\Phi_i, \Phi_j)_{\nabla} := \int_{\U} (\nabla \Phi_i \cdot \nabla \Phi_j) (z) d^2 z.$$
When one explores partially the CLE boundaries, then the effect of the choice of $(\sigma_1, \ldots, \sigma_n)$ will only affect the GFF in $U$ and not in the other (already explored) connected components (that are surrounded by fully explored loops). Furthermore, in $U$, the GFF will then be the sum of a Dirichlet GFF in $U$ with the harmonic function $ \sum_{i}  \sigma_i \Phi_i {\circ \varphi}$. But the GFF is (heuristically) the Boltzmann measure associated to the Dirichlet energy, and the energy of the sum of a field with a harmonic function is just the sum of the Dirichlet energies of the field and of the harmonic function. As a consequence, one can infer that the conditional probability of $(\sigma_1, \ldots, \sigma_n) = (a_1, \ldots, a_n)$ should be proportional  to  $\exp \big(   \sum_{i < j} a_i a_j (\Phi_i, \Phi_j)_{\nabla}/2 \big) = \exp \big(  \sum_{i < j} a_i a_j m_{i,j} \big)$.

Before giving the proof, let us make a brief comment about an idea that does not work so well, namely to try to prove this result inductively over $n$. Indeed,  if one conditions on the first $n-1$ explorations, then the $n$-th exploration itself will actually be strongly correlated with $\sigma_1,\ldots,\sigma_{n-1}$ (for example, if all the $\sigma_i$ agree, then the conditional probability that the $n$-th exploration will get close to all the other explorations will be bigger than if some of the $\sigma_i$ do disagree), so that this approach creates unnecessary difficulties.

\begin {proof}
We will make use of the Cameron-Martin type absolute continuity relations for local sets of the GFF, as in \cite{MR3101840,MR3477777,MR3936643}.
Suppose first that we are given a bounded harmonic function $\Psi$ in a horizontal/vertical rectangle $R$, that has boundary values $0$ except on the right-hand side of the rectangle. We then consider a GFF $h$ with Dirichlet boundary conditions in $R$.
When $\Psi$ is not the zero function, then the law $h + \Psi$ is not absolutely continuous with respect the law of $h$, but this absolute continuity holds when one restricts the fields to the part of $R$ that lies at distance more than some given value from the right-hand side of $R$.  This feature has been exploited in numerous instances in order to make sense of level lines of GFF's with non-zero boundary conditions (see for instance {\cite{MR3101840, MR3477777,MR3936643}}).

One way to get an explicit expression for this Radon-Nikodym derivative it to consider
 some simple boundary-to-boundary curve $S$ in $R$ that separates the right-hand side boundary from some point on its left-hand side boundary. The simple idea is to use some function $\tilde \Psi$ (instead of $\Psi$) that is equal to $\Psi$ on the left-hand side of $S$, but so that the law of (all of) $h + \tilde \Psi$ is absolutely continuous with respect to that of $h$. In order to exploit also  the Markov property of the GFF $h$, it is natural to choose  $\tilde \Psi$ to be the continuous function that is equal to $\Psi$ on $S$, to $0$ on $\partial R$ and that is harmonic in each of the two connected components $R_-$ and $R_+$ of $R \setminus S$. Note that $\tilde \Psi$ is the indeed equal to $\Psi$ in the ``left-hand side'' component $R_-$ of $R \setminus S$, and that $\tilde \Psi$ is still in the Cameron-Martin space of the Dirichlet GFF as it has a finite Dirichlet energy.
By the usual Cameron-Martin result for Gaussian measures (see \cite {MR3936643} for a similar use of it), under the probability measure $\tilde \Pb$ defined to be our original probability measure reweighted by (a constant multiple of) $\exp ( ( h, \tilde \Psi)_\nabla /2)$, the GFF $h$ has the law of the sum of a Dirichlet GFF with $\tilde \Psi$. In particular, its restriction to $R_-$ has the same law as the sum of a Dirichlet GFF with $\Psi$.

We now suppose that $A$ is a local set obtained by exploring the CLE$_4$ of the GFF $h$ starting from a point on the left-hand side $\overline R_-$ of the rectangle along some curve $L$ and stopped when exiting some domain $O$ contained in $\overline R_-$.
This set $A$ is a  thin local set of bounded type  of $h$ (see \cite {MR3936643} for a definition) and one can decompose $h$ as $h_A + h^A$, where $h^A$ is a Dirichlet GFF in the complement of $A$ (i.e., in each of the connected components of the complement of $A$), and $h_A$ is a harmonic function in $R \setminus A$ that takes the values $\pm 2 \lambda$ in each of the connected components fully surrounded by a loop, and the value $\sigma u_A$ in the ``to-be-explored component'' $R_A$ that contains $R_+$, where $u_A$ is the harmonic function in $R_A$ with boundary values $2 \lambda \one_{\xi}$ with $\xi$ being the ``inside'' part of the currently-explored loop.
The CLE$_4$ coupling and the choice of exploration shows that conditionally on the CLE$_4$ exploration, the signs of $h_A$ are independent in each of the connected components of $R \setminus A$. In particular, $\sigma$ is $+1$ or $-1$ with probability $1/2$. So, we define ${\mathcal F}$ to be the $\sigma$-field generated by the exploration of $h$ along $A$ except for the value of $\sigma$, then clearly $ \Pb [ \sigma = +1 | {\mathcal F} ] = \Pb [ \sigma = -1 | {\mathcal F} ] =1/2$.

Then (see  \cite[Proposition 13]{MR3936643}, this argument also appeared in  \cite[Pages 611 -- 613]{MR3477777}),  under $\tilde \Pb$, the set $A$ can be viewed {as} a thin local set of the Dirichlet GFF $\tilde h := h-\tilde{\Psi}$, where the corresponding harmonic function $\tilde{h}_A$
in the complement of $A$ has the same boundary conditions as $h_A -\tilde \Psi$ on the boundary of $R \setminus A$. We therefore swiftly deduce that
$$ \tilde{\Pb} [ \sigma = +1 | {\mathcal F} ] = \exp (( u_A, \tilde \Psi)_{\nabla} / 2 ) \times Y $$
and
$$ \tilde{\Pb} [ \sigma = -1 | {\mathcal F} ] = \exp  (( -u_A, \tilde \Psi)_{\nabla} / 2 ) \times Y $$
almost surely,
where $Y$ is the same ${\mathcal F}$-measurable explicit quantity (involving conditional expectations, $h^A$ and $h_A$ in the other connected components than $R_A$) in both cases.
In particular, we get that almost surely,
$$
   \frac{\tilde{\Pb}[\sigma = +1  |  {\mathcal F}]}{\tilde{\Pb}[\sigma  = -1 | {\mathcal F} ]} =
   \exp ( ( u_A, \tilde \Psi)_\nabla ).
$$
Finally, we can note that
$$ ( u_A, \tilde \Psi)_{\nabla} = ( u_A, \Psi)_{\nabla} $$
which follows directly from the setup (the right-hand side can also be viewed as the limit of the left-hand one when $S$ gets closer and closer to the right-boundary of $R$).

Let us now go back to the setting of the proposition. We choose to first discover the $n-1$ explorations $K_1, \ldots, K_{n-1}$.
Let $U_0:=\U\setminus \cup_{1\le k \le n-1} K_k$. We will rely on the following facts:
\begin{itemize}
\item Conditionally on $K_1, \ldots, K_{n-1}$ and on the labels $\sigma_1, \ldots, \sigma_{n-1}$, $h$ restricted to  $U_0$ is equal to $h_0+ \Psi_0$, where $h_0$ is a GFF  in $U_0$ with zero boundary conditions and   $\Psi_0$ is a harmonic function in $U_0$ with boundary conditions $2 \sigma_1 \lambda, \ldots, 2 \sigma_{n-1} \lambda$ on $\varphi^{-1}(\partial_1), \ldots, \varphi^{-1}(\partial_{n-1})$.  
\item Conditionally on $K_1, \ldots, K_{n-1}$ and on the labels $\sigma_1, \ldots, \sigma_{n-1}$, $K_n$ is a local set of $h_0+ \Psi_0$. Conditionally on $K_n$, the law of the  {GFF $h_0$} in the component with all the explorations on its boundary with boundary conditions $2\sigma_n \lambda$ on $\varphi^{-1}(\partial_n)$ and $0$ elsewhere.
\item If we just condition on $K_1, \ldots, K_{n-1}$ and define the local set $\tilde K_n$ for $h_0$ instead of $h_0 + \Psi_0$ (mind that these two are absolutely continuous when restricted to a domain that $K_n$ is almost surely part of), then the sign $\sigma_n$ of the harmonic function on $\varphi^{-1} (\partial_n)$ is equal to $1$ or $-1$ with conditional probability $1/2$ by  symmetry.
\end{itemize}
We now map $U_0$ conformally onto some rectangle $R$, so that $\varphi^{-1}(\partial_1), \ldots, \varphi^{-1}(\partial_{n-1})$ all belong to the right side of $R$ and  $x_n$ belongs to the left side of $R$. The previous result  provides the conditional probability that the $n$-th label is $+1$, conditionally on $K_1, \ldots, K_{n-1}$ and on the labels $\sigma_1, \ldots, \sigma_{n-1}$.
Indeed,
 that for any given $(a_1, \ldots, a_{n-1})$, 
$$
\frac { p ( a_1, \ldots, a_{n-1} , 1 ) }{ p (a_1, \ldots, a_{n-1}, -1 )}  = \exp \Big(   \sum_{1 \le j \le n-1} 2{a_j} m_{j,n} \Big).
$$
We now use the fact that the role of the $n$ explorations is symmetric, and that we can also apply this reasoning to the case where one instead discovers the $i$-th exploration after the $n-1$ other ones.  We get that for all $i \in \{1, \ldots, n\}$, for all $a_1, \ldots, a_{i-1}, a_i, \ldots, a_n $,
$$ \frac {p (a_1, \ldots, a_{i-1}, 1, a_{i+1}, \ldots, a_n)}{p (a_1, \ldots, a_{i-1}, -1, a_{i+1}, \ldots, a_n)} = \exp \Big(
\sum_{j  \in \{1, \ldots, n\} \setminus \{i\}} 2 {a_j} m_{i,j} \Big). $$
Using this for some $i$ to successively switch some values from $-1$ to $1$ (i.e., those for which $b_j = 1$), we see that this determines in fact
$ p(b_1, \ldots, b_n) / p(-1, \ldots, -1)$ for all $b_1, \ldots, b_n$ and it is easy to check that
the outcome is the formula in Proposition \ref {mainprop2}.
\end {proof}

\subsection {Some facts related to Dynkin's Theorem}
\label {S3.3}

Let us first recall and collect a few facts related to Dynkin's Theorem about squares of the Gaussian Free Field and occupation time of boundary-to-boundary excursions.
The actual version of isomorphism \eqref{eq:dynkin} that we will use can be found in  Sznitman \cite{MR3167123}.

Let us consider a standard GFF $h$ in the unit disc $\U$. When $k$ is a bounded measurable function in $\U$, we will denote by $h(k)$ the value of the GFF tested against $k$.

\subsubsection {The square of the GFF}
To define the (renormalized) square of $h$, one can proceed as follows {(see for instance \cite{MR0489552})}. One first considers a family of smooth mollifiers $\phi_\eps$ that converge to the Dirac mass at $0$ as $\eps \to 0$. In this way, $h_\eps := h * \phi_\eps$ is now a continuous function on $\U$, and it is possible to define its square. It turns out that as $\eps \to 0$,
$$
  \llbracket h^2 \rrbracket  = \lim_{\eps \to 0}\, (h_\eps)^2 -
  \Eb [(h_\eps)^2] $$
exists in $L^2$ when tested against any bounded test function, and that the limit does not depend on the choice of mollifiers $\phi_\eps$. This defines the renormalised square of the GFF. Note that this  process $ ( \llbracket h^2 \rrbracket (k) ,  k \in \mathcal{K})$ is then indexed by the set ${\mathcal K}$ of bounded measurable functions in $\U$, and that all the random variables $ \llbracket h^2 \rrbracket (k))$ are in $L^2$ and centered.

If $\Phi$ denotes a bounded non-identically zero function in $\U$, we can similarly define the renormalized square of $ h + \Phi$ just replacing $h$ by $h+ \Phi$ in the above expression. For our purposes, it will be more suitable to define the process $ \llbracket (h+\Phi)^2 \rrbracket$ as the limit when $\eps \to 0$ of
$$ (h+\Phi)_\eps^2 - \Eb [ (h_\eps)^2]$$
(observing that this is then no longer a centered process).
Noting simply that
$$((h+ \Phi)_\eps)^2 = (h_\eps + \Phi_\eps)^2 =  h_\eps^2 + \Phi_\eps^2 + 2 h_\eps \Phi_\eps , $$ we get that
$$ \llbracket (h + \Phi)^2 \rrbracket = \llbracket h^2 \rrbracket + 2 h \Phi + \Phi^2 .$$
In other words, for each $k$ in ${\mathcal K}$,
$$ \llbracket (h+\Phi)^2 \rrbracket (k) =
\llbracket h^2 \rrbracket (k) + 2 h ( k \Phi) + \int_\U \Phi^2 (z) k(z) d^2 z. $$

\subsubsection {Weighted GFF}
Fix a smooth compactly supported non-negative function $k$ defined in $\U$.
When $h$ is a standard GFF  in $\U$ (i.e., with generator $-\Delta/2$), one can reweight its law by $$\exp(-\llbracket h^2\rrbracket (k) / 2 ),$$ and it is then easy to see (see for instance {\cite[Proposition 4.3] {Dub}}) that under this reweighted
measure, the field $h$ is still Gaussian, and can be viewed as a Gaussian Free Field associated to the Brownian motion killed at rate $k$, i.e., with generator  $-\Delta/2 + k$.

Let us denote this field by $h_k$ (not to be confused with our previous notation $h_\eps$ -- but we will not use it anymore in the rest of the paper), and the corresponding Green's function by $G_k$.
We see  that
\begin {eqnarray*}
 \lefteqn { 
\Eb \left[ \exp \left( -  \frac 1 2 \llbracket (h+\Phi)^2 \rrbracket (k) \right) \right]}\\
&&= \exp \left( -  \frac 1 2 \int_\U \Phi^2 (z) k(z) d^2 z \right) \times
\Eb  \left[\exp \left( - \frac 1 2  \llbracket h^2 \rrbracket (k) \right) \right] \times \Eb [ \exp (- h_k ( k \Phi )) ],
\end {eqnarray*}
and we can note that (by looking at the variance of the centered Gaussian variable $h_k ( k \Phi)$), the final term has the explicit form
$$
 \Eb [ \exp (-h_k ( k \Phi )) ] = \exp \left( \frac 1 2 \int_{\U \times \U} G_k (w, z) k(w) k(z) \Phi(w) \Phi (z)  d^2 z d^2 w \right).$$

 \subsubsection {Using Dynkin's isomorphism}
To state Dynkin's Theorem (i.e., one version of it), we will use here the standard Brownian excursion measure in $\U$ (that is supported on the set of Brownian paths that start and end on $\partial \U$). When $\partial$ is a closed subset of $\partial \U$, we denote $\mu_\partial$ to be the excursion measure $\mu$ restricted to the excursions that have both endpoints in $\partial$, and we let $\Phi_\partial$ denote the harmonic extension in $\UU$ of {$2 \lambda \one_\partial$} defined on $\partial \U$.
We then let $T_\beta^\partial$ be the
occupation field of a Poisson point process with intensity $\beta
\mu_\partial$ that is independent of the GFF $h$.
Then, Dynkin's Theorem (see for instance \cite{MR3167123}) states that
\begin{align}\label{eq:dynkin}
    \frac{1}{2}\llbracket (h+u\Phi_\partial)^2 \rrbracket \stackrel{d}{=}
    \frac{1}{2} \llbracket h^2 \rrbracket + T^\partial_\beta
\end{align}
where {$\beta = u^2 / 4$.   Let us also observe that $u^2=1 $ corresponds to the value $\beta = 1/4$ -- these are actually the values that we will be working with. }
In particular, combining this with the considerations and results of the previous section shows that for any bounded non-negative function $k$ in $\U$,
\begin {eqnarray*}
\lefteqn {\Eb [ \exp ( - T_\beta^\partial (k) )]} \\
& =& \frac {1}{\Eb [ \exp (-  \llbracket h^2 \rrbracket (k)  /2)]}\times
 \Eb \left[ \exp \left( - \frac 1 2 \llbracket h^2 \rrbracket (k) - u h ( \Phi_\partial k) \right) \right] \\
 && \times \exp \left( - \frac {u^2} 2 \int_\U \Phi_\partial^2(z)  k(z) d^2 z \right) \\
 &=&
 \Eb [ \exp (  - u h_k( \Phi_\partial k) )] \times \exp \left( - \frac {u^2} 2 \int_\U \Phi_\partial^2(z)  k(z) d^2 z \right) \\
&=&
  \exp \left( - \frac {u^2} 2 \int_\U  \Phi_\partial^2 (x)  k (x) d^2x \right)
{\times}   \exp \left(   \frac {u^2} 2 \int_{\U \times \U}  G_k(w,z)\Phi_\partial (w) k(w) \Phi_\partial (z) k(z) \,d^2 w d^2 z \right)
\end {eqnarray*}

\subsubsection {Regrouping terms}
Finally, noting that
$$
\Eb [ \exp ( - T_\beta^{\partial_i \cup \partial_j} (k) )]
=
\Eb [ \exp ( - T_\beta^{\partial_i} (k) )]
\Eb [ \exp ( - T_\beta^{\partial_i} (k) )]
\Eb [ \exp ( - T_\beta^{\partial_i \leftrightarrow \partial_j} (k) )],
$$
where $T_\beta^{\partial_i \leftrightarrow \partial_j}$ corresponds to a Poisson point process of Brownian excursions with one endpoint in $\partial_i$ and one endpoint in $\partial_j$,
we get that
\begin {eqnarray*}
\lefteqn {   \Eb[\exp ( -T^{\partial_i \leftrightarrow \partial_j}_\beta (k))]}\\
&=&
  \exp \left( - {u^2} \int_\U  \Phi_{i}(x) \Phi_{j}(x)  k (x) d^2x +   u^2 \int_{\U \times \U}  G_k(w,z)\Phi_{i} (w) k(w) \Phi_{j} (z) k(z) \,d^2w d^2 z \right).
\end {eqnarray*}
We can note that (for instance by applying this to $yk$ and letting $y \to 0$) this implies that
$$ {u^2} \int_\U  \Phi_{i}(x) \Phi_{j}(x)  k (x) d^2x = \Eb [ T^{\partial_i \leftrightarrow \partial_j}_\beta (k) ].$$

An observation that will be handy later on is to note that by the standard Laplace transform for Poisson point processes,
$$ \Eb [ \exp ( - T_\beta^{\partial_i \leftrightarrow \partial_j} (k) )]
= \exp ( -  \beta \mu_{i,j} ( 1 - \exp ( -T_e (k) ))$$
(where $T_e$ is now the occupation time of the excursion $e$). 
Hence (recalling that the total mass $|\mu_{i,j}|$ of $\mu_{i,j}$ is finite),
we get that $  {\beta \mu_{i,j} ( \exp ( -T_e (k) ))}$ is equal to
\begin {equation}
 \beta |\mu_{i,j}| - {u^2} \int_\U  \Phi_{i}(x) \Phi_{j}(x)  k (x) d^2x
+  u^2 \int_{\U \times \U}  G_k(w,z)\Phi_{i} (w) k(w) \Phi_{j} (z) k(z) \,d^2w d^2 z .
\label {mu_ij}
\end {equation}

\subsection {Parity-constrained excursions}
\label {pce}
\label {S3.4}

We are now almost ready to conclude. We will want to compare the law of occupation times that we have just obtained for our set of excursions with that of parity constrained Poisson point processes. Let us first describe the latter in this short section:
Let us now consider a Poisson point process of excursions with intensity $\beta \mu$ restricted to the set of excursions such that one endpoint is in some $\partial_i$ and the other endpoint is in another $\partial_j$ for some $j \in \{ 1, \ldots, n \} \setminus \{ i \}$.
We let $\n_{i,j}$ denote the number of excursions joining $\partial_i$ and $\partial_j$, write $T_{i,j}$ for its occupation field and we then denote by $\n_i$ the number of excursions that have one endpoint in $\partial_i$. In other words, $\n_i := \sum_{j \not=i } \n_{i,j}$.
We finally will condition this Poisson point process of excursions on the event
$$ {\mathcal E} := \{ \n_i \hbox { is even for all  } i \}.$$

Before that, we can note that $\n_{i,j}$'s are independent Poisson random variables with respective means $ \beta |\mu_{i,j}|$, and that for each individual excursion is then chosen independently according to $\mu_{i,j} / | \mu_{i,j}|$, so that
$$
\Eb [ \exp ( - T_{i,j} (k) )]
= \sum_{n_{i,j} \ge 0} e^{-|\beta \mu_{i,j}|} \frac { \bigl( \beta \mu_{i,j} ( \exp ( - T_e(k)))\bigr)^{n_{i,j}}}{n_{i,j}!},
$$
where $T_e (k) = \int_0^\tau k(e(s)) ds$ for an excursion $e$ with time-length $\tau$.
We can use the standard trick (similar to the random current representation of the Ising model) that
\begin {eqnarray*}
\lefteqn{\sum_{(a_1, \ldots, a_n) \in \{ -1, 1\}^n }
\prod_{i< j}
\sum_{n_{i,j} \ge 0} e^{-|\beta \mu_{i,j}|} \frac {\bigl(a_i a_j \beta \mu_{i,j} ( \exp ( - T_e(k)))\bigr)^{n_{i,j}}}{n_{i,j}!}}
\\
&=&
{\sum_{(a_1, \ldots, a_n) \in \{ -1, 1\}^n }
\sum_{(n_{i,j}) \ge 0} \left(\prod_i a_i^{n_i} \right) \prod_{i<j}   e^{-|\beta \mu_{i,j}|}  \frac {\bigl( \beta \mu_{i,j} ( \exp ( - T_e(k)))\bigr)^{n_{i,j}}}{n_{i,j}!}}
\\
& =& {2^n}\, \Eb \left[ \one_{\mathcal E} \prod_{i < j} \exp ( -T_{i,j} (k)) \right]
.\end {eqnarray*}
Hence, we can readily conclude that
$$ \Eb \left[ \prod_{i < j} \exp ( -T_{i,j} (k))  | {\mathcal E} \right]
= \frac {1}{Z'}
\sum_{a_1, \ldots, a_n} \prod_{i< j} \exp \bigl( \beta  a_i a_j \mu_{i,j} (\exp (- T_e (k))) \bigr)$$
where the constant $Z'$ is chosen so that the quantity is $1$ for $k=0$.

\subsection {Wrapping up the GFF computation}
\label {wrapping}
\label {S3.5}
On the GFF side, we can do the computation of the occupation time. More precisely, we can compute the conditional expectation of
$\exp ( - \llbracket h^2 \rrbracket (k))$ given the CLE exploration, when $h$ is coupled to the loop-soup. We will restrict ourselves to the bounded functions $k$  on $\U$.

We can then decompose according to values of $\sigma_1, \ldots, \sigma_n$, and we obtain the following expression for this conditional expectation:
\begin {eqnarray*}
\lefteqn { \sum_{(a_i) } p(a_1, \ldots, a_n) \times \Eb [ \exp(-  \llbracket (h+  \sum a_i \Phi_i)^2  \rrbracket (k)  )]}
\\
&= & \sum_{(a_i)}
\frac {1}{Z}
\exp \left(  \sum_{i<j} a_i a_j m_{i,j} \right) \times
\Eb \left[  \exp (- \llbracket h^2 \rrbracket  (k)) \right] \\
&&\times
\exp \left( -  \frac{1}{2} \int_\U \Phi_{(a_i)}^2 (z) k(z) d^2 z \right) {\times}
\exp \left( \frac{1}{2} \int_{\U \times \U} G_k (w, z) k(w) k(z) \Phi_{(a_i)}(w) \Phi_{(a_i)} (z)  d^2 z d^2 w \right)
\end {eqnarray*}
where $\Phi_{(a_i)}  = \sum_{i} a_i \Phi_i$ and $h$ is a GFF in $\U$.
By expanding the products and regrouping terms  {(for each fixed $a_1, \ldots, a_n$),} we get the expression
$$ \sum_{(a_i)} [ {\frac{1}{Z} \times (L) \times \prod_{1\le i \le n} (L)_i \times \prod_{1 \le i<j\le n} (L)_{i,j}}]$$
where
\begin {eqnarray*}
{(L)} &:=& \Eb [ \exp (- \llbracket h^2 \rrbracket  (k) ) ],
\\
{(L)_i}  &:=&
\exp \left( -   \frac{1}{2} \int_\U \Phi_i^2 (z) k(z) d^2 z \right)\\
&& {\times}
\exp \left(  \frac{1}{2} \int_{\U \times \U} G_k (w, z) k(w) k(z) \Phi_i(w) \Phi_i (z)  d^2 z d^2 w \right)
\\
{(L)_{i,j}} &:= & \exp (  a_i a_j m_{i,j} ) {\times}
\exp \left(  a_i a_j \int_\U \Phi_{i}(z) \Phi_j (z) k(z) d^2 z \right) \\
&& {\times}
\exp \left( a_i a_j  \int_{\U \times \U} G_k (w, z) k(w) k(z) \Phi_i(w) \Phi_j (z)  d^2 z d^2 w \right).
\end {eqnarray*}
We recognize the first one as giving the renormalized occupation time of the Brownian loop-soup, and the second one as the (non-renormalized) occupation time of a Poisson point process of excursions from $\partial_i$ to $\partial_i$.

Using (\ref {mu_ij}) for {$u=1$}, we see that
$$ {(L)_{i,j}}
=  \exp \bigl(  a_i a_j  \beta \mu_{i,j} (\exp (- T_e (k))) \bigr)$$
for $\beta = 1 / 4$,
so that indeed, $\prod_{i<j} {(L)_{i,j}}$ (summed over all configurations $(a_1,\ldots,a_n) \in \{-1,1\}^n$) corresponds to the occupation time of the parity-conditioned Poisson point process of excursions (with end-points in different $\partial_i$'s).

Recall that we know that conditionally on the explorations, the not-yet discovered loops in $\U$ do form a Brownian loop-soup that is independent of the collections of loops intersecting  $\partial \U$.  We can therefore conclude that conditionally on the collection of loops not intersecting $\partial$, the occupation time field of the loops that do intersect $\partial$ has indeed the law of a parity-conditioned Poisson point process of excursions.

\subsection {Concluding the proof using \cite {MR3896865}}
\label {S3.6}
The result \cite[Proposition 1.1]{MR3896865} states that if one knows the law of the trace of the union of the excursions in a {locally finite} point process (i.e., not necessarily a Poisson point process) of Brownian excursions in a domain, then one knows the law of this point process.

In the present case, we have shown that the (conditional, given the explorations) law of the occupation times of the union of the excursions of the boundary touching excursions in $\U$ is the same as that of the (occupation times of the union of the excursions in the) parity-conditioned Poisson point process of Brownian excursions.  {Note that in both cases, the occupation times determine the trace (because the intersection of these union of excursions with any disk $\lambda {\overline{\U}}$ for $\lambda < 1$ is compact), so that one gets also the identity between the traces of these point processes.}

To apply the result of \cite {MR3896865}, it therefore remains to check that these excursions are indeed Brownian excursions (i.e., when conditioned on their endpoints, the law of each excursion is that of {an independent} Brownian excursion with these endpoints) {and that there are finitely many of them which have diameter at least $\delta$, for any $\delta>0$}. This can be shown as in {\cite[Lemma 9] {MR3896865}} (where the excursions inside a ``fully discovered boundary'' were studied) using the resampling property of the Brownian loop-soup. This concludes the proof of Theorem \ref {mainthm}.

\subsection {More general explorations}
\label {MoreGeneral}
We now explain how to generalize Theorem \ref {mainthm} to more general explorations. We will only outline one possible way to proceed, and leave the details to the interested reader.
Our choice of the $n$ deterministic lines was rather arbitrary, as well as the orientation choice of each loop. A first remark is that the proof works in exactly the same way if one replaces the deterministic collections $L_1, O_1, \ldots, L_n, O_n$ by random ones that are independent from each other and from the loop-soup and remain disjoint. Similarly, one can also toss a new coin at each time the exploration hits a loop to decide which way to trace it. More generally, one can use any set of $n$ Markovian explorations of the CLE$_4$, as long as the ``rules'' that determine these $n$ explorations (and when to stop them) are independent from each other -- this includes using level lines of a GFF coupled to the CLE$_4$ (see Remark \ref {levellines} below).

One useful further observation is that Theorem \ref {mainthm} shows that the conditional distribution
of the loop-soup configuration (i.e., the excursions plus the remaining loops) after the $n$ explorations is a conformally invariant function of the domain $U$ with the $n$ boundary arcs (the excursions form a parity-constrained point process and the other loops form an independent loop-soup), where the roles of the boundary arcs is in fact symmetric.

As a consequence, we see that if one actually starts with such a configuration  (with $n$ boundary arcs, a loop-soup and a parity constrained collection of excursions) and continues exploring the last strand for a little while (say along $\tilde L_n$ that extends $L_n$), the conditional law of the configuration after this last additional exploration is again that of a loop-soup with a parity constrained collection of excursions. We can note that this will hold for any given additional $\tilde L_n$. Hence, this would also hold if this $\tilde L_n$ was chosen using the information provided by the explorations along $L_1, \ldots, L_n$.

Combining all this allows to in fact ``successively switch'' from one exploration strand to the other, and to use additional randomness (as long as it is independent of the ``to be discovered configuration'') in order to decide how/where to explore next. In this way, one can then essentially approximate a large-class of reasonable ``Markovian'' exploration from $n$ marked points, and then obtain the same conditional distribution for the remaining configuration.

This provides also a way, for each given $x$, to define iteratively two explorations stopped in such a way (by first choosing the first one, and then to stop the second one at a time chosen depending on the first one) that $\partial_1$ and $\partial_2$ would be a conformal rectangle with aspect ratio $x$ with positive probability.

\begin {remark}
\label {levellines}
Recall that one natural, conformally invariant and special way to explore a CLE$_4$ is along SLE$_4$-type curves (see \cite{MR2494457}) which can be interpreted as level-lines of a GFF coupled to this CLE$_4$ (see \cite{MR3936643,MR3708206}).  Recall that if one gives i.i.d.\ orientations to each loop in the CLE$_4$, with probability $1/2$ for each orientation, then there is a continuous Markovian exploration of the \emph{oriented} CLE$_4$ along a random curve $\gamma$ started from any boundary point $x$, so that we again trace each loop that the curve encounters in the given orientation, and in the order that $\gamma$ encounters them.
However, if we would like to use such an exploration for the $n\ge 2$ case, then we need to use different independent coin tosses for the loops in the different explorations (this was also explained in \cite[Section 3]{MW}), even if they actually do turn out to correspond to the same loop. In other words, each exploration uses a different (conditionally independent given the CLE) GFF.
\end {remark}

\subsection {The parity of rectangle-crossings}
\label {Sfinal}

Let us now explain finally explain how to use Theorem \ref {mainthm} to deduce the results about rectangle crossing mentioned in the introduction,

As a warm-up, let us first indicate why the coupling of the CLE$_4$, the GFF and the loop-soup indicate that the parity of the number of crossings is not a deterministic function of the occupation field.
Let us first perform an exploration from two points of the CLE$_4$, as in Theorem \ref {mainthm}. We call ${\mathcal H}$ the corresponding $\sigma$-field. The obtained set $U$ with the two boundary arcs can be conformally mapped onto  a rectangle $R$ in such a way that the ``inner sides'' of $\xi_1$ and $\xi_2$ are respectively mapped onto the two vertical sides $V_1$ and $V_2$ of $R$. We call $X$ the (random) aspect ratio of this rectangle.

Theorem \ref {mainthm} then describes the remaining to be discovered part of the loops as the (conformal preimage) of a critical loop-soup $\Gamma$  in $R$ together with an independent parity-constrained collection $\Sigma$ of Brownian excursions from $V_1$ and $V_2$ as defined by Theorem~\ref{mainthm}.
In the natural coupling of the loop-soup with the GFF, the renormalised occupation time of the union of these two will correspond to the renormalised square of the GFF.

We also know the probability $c(X)$ that the two arcs $\xi_1$ and $\xi_2$ will be in the same CLE$_4$ loop or not (this is for instance part of \cite {MW} in the special $\kappa=4$ case, as described in Section 3 of that paper). We let $A$ denote the event that they are in the same loop. Then, we know that
\begin {itemize}
 \item When $A$ holds, then in coupling of the CLE$_4$ with the GFF, both arcs necessarily have the same sign (i.e., $\sigma_1 = \sigma_2$).
 \item When $A$ does not hold, then the signs of the two arcs will be independent, and have probability $1/2$ to be the same.
\end {itemize}
So, we can conclude that the law of the renormalised square of the GFF conditionally on $\{ \sigma_1 = \sigma_2 \}$ and ${\mathcal H}$ is absolutely continuous with respect to the unconditional law of the renormalised square of the GFF given ${\mathcal H}$.

We now define ${\mathcal H}'$ to be the $\sigma$-field generated by ${\mathcal H}$ and $(\sigma_1, \sigma_2)$. This corresponds to a GFF exploration. In particular, conditionally on ${\mathcal H}'$, the law of the GFF in $U$ is just that of a GFF with boundary conditions $2\lambda (\sigma_1 \one_{V_1} + \sigma_2 \one_{V_2})$. By Dynkin's isomorphism theorem, when $\sigma_1 = \sigma_2$, the law of its renormalized square is that of the sum of the renormalized occupation time of a loop soup with the occupation time of a Poisson point process of excursions away from $V_1 \cup V_2$ (with no parity constraint). We can therefore finally conclude that for almost all $x$ (with respect to the law of $X$), the law of the sum of a renormalised square of a GFF with the occupation time of an independent  Poisson point process of excursions from $V_1 \cup V_2$ in $R$ with no parity constraint is absolutely continuous with {respect to} the law of the sum of  a renormalised square of a GFF with the occupation time of an independent  Poisson point process of excursions from $V_1 \cup V_2$ in $R$ with parity constraint.

To deduce that this property holds in fact for all positive $x$, one can for instance first show using density properties for the loop-soup that the possible laws of $X$ (letting the choices of explorations vary -- in the spirit of Appendix B of \cite {MR4055986}) imply that the result holds for a dense set of $x$ in $(0,\infty)$ (with a Radon-Nikodym derivative that is bounded and bounded locally from below locally with respect to $x$ in that set), and then by approximation (letting $x_n \to x$ for any given value of $x$ and $x_n$ such that the result holds) conclude that this holds as well for all $x$. We leave details to the reader. Another approach is to use the generalizations to iterated Markovian explorations as described in the previous section.

\medbreak
Let us now explain that a stronger result actually holds. We still work with the rectangle $R$ with aspect ratio $x$, and consider a Poisson point process of left-right Brownian excursion (with the same particular intensity) and an independent loop-soup. We denote by $A$ the event that the left-side and the right-side of the rectangle are connected by a cluster of the union of the loop-soup and of the excursions. We let $m$ denote the mass of Brownian excursions from the left and to the right side of the rectangle (normalized so that the number $N$ of excursions is a Poisson random variable with mean $m$). We denote by $E$ (respectively $O$) the event that $N$ is even (resp. odd).  Clearly, the complement of $A$ can occur only when $N=0$, and we also obviously have that
$$ \Pb [ E ] =  \frac {1 + \exp (-2m)}{2},  \Pb [ O] = \frac{ 1- \exp (-2m)}{2}  \hbox { and } \Pb [ E ] - \Pb [ O] = \exp (-2m).$$
One can also compute the value of $\Pb [ E \cap A^c] = \Pb [ \{ N =0 \} \cap A^c]  $ simply and directly via SLE$_4 ( \rho)$ and restriction measure considerations (as for instance in \cite {MR3941462}), but we also have another simple way to proceed here:
 Indeed, by Theorem \ref{mainthm}, we see that the probability that $\sigma_1 \not= \sigma_2$ in Proposition \ref {mainprop2} equals the probability of $A^c$ given $E$ multiplied by $1/2$ (since conditionally on the two sides being in different clusters, we are sampling two independent random signs for the two clusters). Rearranging and using the proposition yields
$$\Pb [ E \setminus A | E ]
= \frac{2 \exp (-m)}{\exp (-m) + \exp (m)} =
\frac {2 \exp (-2m) }{   1 + \exp(-2m)} .$$
Hence, 
$$ \Pb [ E \cap A ] = P  [E]  (1  - \Pb [E \setminus A | E ])   = \frac{1- \exp( -2m)}{2} , $$
which happens to be equal to $\Pb [ O]$.
In other words: 
\begin {lemma}
Given that the left and the right side are joined by a cluster, the conditional probability that the number of left-right crossing excursions is even is equal to $1/2$ (and therefore the conditional probability that this number is odd is equal $1/2$ as well).
\end {lemma}

This suggests an occupation-measure preserving bijective switching mechanism from $E \cap A$ onto $O$. Indeed:
\begin {theorem}
The laws of the renormalized occupation time measures of the union of the excursions and the loop-soup, when conditioned on $E \cap A$ and when conditioned on $O$ are equal.
\end {theorem}

\begin {proof}
The proof reuses similar ideas as that of Theorem \ref {mainthm}. Since $\Pb [E \cap A] = \Pb [ O]$ and by Theorem \ref {mainthm}, we see that it is sufficient to check that for any function $k \in {\mathcal K}$, the conditional expectation of  $\one_{F} \exp ( - [[(h+ \sigma_1 \Phi_1 + \sigma_2 \Phi_2)^2]] (k) )$ (when $F$ denotes the event that the two partially explored arcs are part of the same CLE cluster) given the exploration is a multiple (where the constant does not depend on $k$)  of the expression obtained for the union of a loop-soup and a Poisson point process of excursions, restricted to the event $O$ where the number of left-right crossing excursions is odd.

The computation for the latter part goes along as in the proof of Theorem \ref {mainthm}, except that the parity-constrained excursions computation from $\partial_1$ to $\partial_2$ in  Section \ref{pce} has to be replaced by the sum over an odd number of excursions, so that the contribution to the Laplace transform of these odd crossing excursions becomes (where $m$ is the total mass of $\beta \mu_{1,2}$ for  {$\beta= 1/ 4$}),
$$ e^{-m} \left( \exp ( \beta \mu_{1,2} ( e^{-T_e (k)} )) - \exp ( -\beta \mu_{1,2} ( e^{-T_e (k)} )) \right) =  2e^{-m} \sinh(\beta \mu_{1,2} ( e^{-T_e (k)} )).$$
Recall that the expression for $\beta \mu_{1,2} ( \exp (- T_e (k)))$ is still given by (\ref {mu_ij}) for $u=1$ (i.e., $\beta = 1/ 4 $). So, altogether, one gets a constant (i.e., independent of $k$) times
\begin {equation}
 \label {finalexpression}
\sinh \left( m -  \int_{\U} \Phi_1 (x) \Phi_2 (x) {k(x)} d^2 x + \int_{\U \times \U} G_k (w,z) \Phi_1 (w) \Phi_2 (z) k(w) k (z) d^2 w d^2 z \right).\end {equation}

For the computation of the Laplace transform on the exploration side, we can use the coupling with the GFF to notice that the contribution of the event that $\sigma_1 = \sigma_2$ {\em  and} that the two sides are not joined by a cluster is identical to the contribution of the event when $\sigma_1 = -\sigma_2$ {(see also Section \ref{Sfinal} where we made the same observation)}. As a consequence, restricting to the event $F$ amounts to just changing the sign of the contribution of $\sigma_1 = - \sigma_2$ in the sum
(formally, we use that $\Omega$ is the disjoint union of  $F = F \cap \{ \sigma_1 = \sigma_2\}$, $F^c \cap \{\sigma_1 = \sigma_2\}$ and  $F^c \cap \{ \sigma_1 \neq \sigma_2\} =
\{\sigma_1 \neq \sigma_2\}$ and observe that  the second and third event have the same probability, so that when changing the sign of the contribution of the last event, only the probability of $F$ remains).
One therefore expands just as in Section \ref {wrapping} for $n=2$ explorations, except that one puts a minus sign in front of the terms $p(1, -1)$ and $p(-1, 1)$.
In the factorization, the terms ${(L)}$, ${(L)_i}$ and ${(L)_{i,j}}$ remain unchanged, but when one sums them over all choices of $a_1, a_2$, one has to sum the terms $a_1a_2 {(L)_{1,2}}$ instead of the terms ${(L)_{1,2}}$. This sum then turns out to indeed become a constant (that does not depend on $k$) times the expression (\ref {finalexpression}).

Given that the two events $A \cap E$ and $O$ have the same probability, the multiplicative constants actually match (this corresponds to the choice $k=0$). This shows the result for almost all aspect-ratio with respect to the law of the aspect-ratio of the conformal rectangle obtained for each given exploration mechanism. To get the result for each fixed aspect-ratio and therefore conclude the proof, we can proceed as outlined in the previous sections.
\end {proof}

This raises the question of whether one can construct explicitly this bijective occupation-time measure preserving bijection from $E \cap A$ onto $O$. One could tentatively imagine {this} in the spirit of the { ideas} appearing in Section \ref {S2}. In some sense,  in the setting of Figure \ref {fig:loop-soup-extra}, this could for instance correspond to only performing the switches ``within the additional excursion'' (and the loop-cluster it meets) but not on its chosen way back on the boundary. It is then entertaining to think of this tentative Markov chain to move from even number of crossings to odd number of crossings and vice-versa by $\pm 1$ at each step, where the marginal measure on the number of crossings of the invariant measure is a multiple of the Poisson distribution of parameter $m$ except at the origin.

\subsubsection*{Acknowledgements}
WQ is partially supported by a GRF grant from the Research Grants Council of the Hong Kong SAR (project CityU11305823). WW is supported by a Royal Society Research Professorship.

\bibliographystyle{plain}

\begin{thebibliography}{99}

\bibitem{MR4574830}
\'{E}lie A\"{\i}d\'{e}kon, Nathana\"{e}l Berestycki, Antoine Jego, and Titus
  Lupu.
\newblock Multiplicative chaos of the {B}rownian loop soup.
\newblock {\em Proc. Lond. Math. Soc. (3)}, 126(4):1254--1393, 2023.

\bibitem{MR4399157}
Juhan Aru, Titus Lupu, and Avelio Sep\'{u}lveda.
\newblock Extremal distance and conformal radius of a {$\rm CLE_4$} loop.
\newblock {\em Ann. Probab.}, 50(2):509--558, 2022.

\bibitem{aru2023excursion}
Juhan Aru, Titus Lupu, and Avelio Sepúlveda.
\newblock Excursion decomposition of the 2d continuum {GFF}, 2023.
\newblock arXiv 2304.03150.

\bibitem{MR3936643}
Juhan Aru, Avelio Sep\'{u}lveda, and Wendelin Werner.
\newblock On bounded-type thin local sets of the two-dimensional {G}aussian
  free field.
\newblock {\em J. Inst. Math. Jussieu}, 18(3):591--618, 2019.

\bibitem{MR2435854}
Vincent Beffara.
\newblock The dimension of the {SLE} curves.
\newblock {\em Ann. Probab.}, 36(4):1421--1452, 2008.

\bibitem{MR1163370}
Leo Breiman.
\newblock {\em Probability}, volume~7 of {\em Classics in Applied Mathematics}.
\newblock Society for Industrial and Applied Mathematics (SIAM), 1992.
\newblock Corrected reprint of the 1968 original.

\bibitem{MR1062056}
Krzysztof Burdzy and Gregory~F. Lawler.
\newblock Nonintersection exponents for {B}rownian paths. {II}.\ {E}stimates
  and applications to a random fractal.
\newblock {\em Ann. Probab.}, 18(3):981--1009, 1990.

\bibitem{Dub}
Julien Dub{\'e}dat.
\newblock S{LE} and the free field: partition functions and couplings.
\newblock {\em J. Amer. Math. Soc.}, 22(4):995--1054, 2009.

\bibitem{Duplantier}
Bertrand Duplantier.
\newblock Conformally invariant fractals and potential theory.
\newblock {\em Phys. Rev. Lett.}, 84(7):1363--1367, 2000.

\bibitem{dynkin1983markov}
Eugene~B. Dynkin.
\newblock Markov processes as a tool in field theory.
\newblock {\em Journal of Functional Analysis}, 50(2):167--187, 1983.

\bibitem{dynkin1984gaussian}
Eugene~B. Dynkin.
\newblock Gaussian and non-gaussian random fields associated with markov
  processes.
\newblock {\em Journal of Functional Analysis}, 55(3):344--376, 1984.

\bibitem{eisenbaum2000ray}
Nathalie Eisenbaum, Haya Kaspi, Michael~B. Marcus, Jay Rosen, and Zhan Shi.
\newblock A {R}ay-{K}night theorem for symmetric {M}arkov processes.
\newblock {\em Ann. Probab.}, 28(4):1781--1796, 2000.

\bibitem{GLQ_loopsoup}
Yifan Gao, Xinyi Li, and Wei Qian.
\newblock Multiple points on the boundaries of {B}rownian loop-soup clusters.
\newblock {\em Ann. Probab.}, to appear.

\bibitem{MR3786302}
Ewain Gwynne, Jason Miller, and Xin Sun.
\newblock Almost sure multifractal spectrum of {S}chramm-{L}oewner evolution.
\newblock {\em Duke Math. J.}, 167(6):1099--1237, 2018.

\bibitem{JLQ_loopsoup}
Antoine Jego, Titus Lupu, and Wei Qian.
\newblock Conformally invariant fields out of {B}rownian loop soups, 2023.
\newblock arXiv 2307.10740.

\bibitem{MR2644878}
Richard Kiefer and Peter M\"orters.
\newblock The {H}ausdorff dimension of the double points on the {B}rownian
  frontier.
\newblock {\em J. Theoret. Probab.}, 23(2):605--623, 2010.

\bibitem{MR1117680}
Gregory~F. Lawler.
\newblock {\em Intersections of random walks}.
\newblock Probability and its Applications. Birkh\"{a}user, 1991.

\bibitem{Law}
Gregory F. Lawler.
Geometric and fractal properties of Brownian motion and random walk paths in two and three
dimensions.
in {\em Random Walks} (Budapest, 1998), 219--258, Bolyai Soc. Math. Stud., 1999.

\bibitem {LSW1}
Gregory~F. Lawler, Oded Schramm, and Wendelin Werner.
Values of Brownian intersection exponents, II: Plane exponents.
\newblock {\em Acta Math.} 187:275--308, 2001.



\bibitem{MR1849257}
Gregory~F. Lawler, Oded Schramm, and Wendelin Werner.
\newblock The dimension of the planar {B}rownian frontier is {$4/3$}.
\newblock {\em Math. Res. Lett.}, 8(4):401--411, 2001.

\bibitem{MR1961197}
Gregory~F. Lawler, Oded Schramm, and Wendelin Werner.
\newblock Analyticity of intersection exponents for planar {B}rownian motion.
\newblock {\em Acta Math.}, 189(2):179--201, 2002.

\bibitem{MR1992830}
Gregory~F. Lawler, Oded Schramm, and Wendelin Werner.
\newblock Conformal restriction: the chordal case.
\newblock {\em J. Amer. Math. Soc.}, 16(4):917--955, 2003.

\bibitem{MR2045953}
Gregory~F. Lawler and Wendelin Werner.
\newblock The {B}rownian loop soup.
\newblock {\em Probab. Theory Related Fields}, 128(4):565--588, 2004.

\bibitem{MR1229519}
Jean-Fran{\c{c}}ois Le~Gall.
\newblock Some properties of planar {B}rownian motion.
\newblock In {\em \'{E}cole d'\'{E}t\'e de {P}robabilit\'es de {S}aint-{F}lour
  {XX}---1990}, volume 1527 of {\em Lecture Notes in Math.}, pages 111--235.
  Springer,  1992.

\bibitem{le2010markov}
Yves Le~Jan.
\newblock Markov loops and renormalization.
\newblock {\em The Annals of Probability}, pages 1280--1319, 2010.

\bibitem{MR2815763}
Yves Le~Jan.
\newblock {\em Markov paths, loops and fields}, volume 2026 of {\em Lecture
  Notes in Mathematics}.
\newblock Springer,  2011.
\newblock [Lectures from the 38th Probability Summer School held in Saint-Flour in
  2008].

\bibitem{MR3502602}
Titus Lupu.
\newblock From loop clusters and random interlacements to the free field.
\newblock {\em Ann. Probab.}, 44(3):2117--2146, 2016.

\bibitem{MR3941462}
Titus Lupu.
\newblock Convergence of the two-dimensional random walk loop-soup clusters to
  {CLE}.
\newblock {\em J. Eur. Math. Soc.}, 21(4):1201--1227, 2019.

\bibitem{marcus2006markov}
Michael~B. Marcus and Jay Rosen.
\newblock {\em Markov processes, Gaussian processes, and local times}.
\newblock Cambridge University Press, 2006.

\bibitem{MS}
Jason Miller and Scott Sheffield.
\newblock {CLE}(4) and the {G}aussian {F}ree {F}ield.
\newblock {\em Private communication}.

\bibitem{MR3477777}
Jason Miller and Scott Sheffield.
\newblock Imaginary geometry {I}: interacting {SLE}s.
\newblock {\em Probab. Theory Related Fields}, 164(3-4):553--705, 2016.

\bibitem{MR3708206}
Jason Miller, Scott Sheffield, and Wendelin Werner.
\newblock C{LE} percolations.
\newblock {\em Forum Math. Pi}, 5:e4, 102, 2017.

\bibitem{MR4055986}
Jason Miller, Scott Sheffield, and Wendelin Werner.
\newblock Non-simple {SLE} curves are not determined by their range.
\newblock {\em J. Eur. Math. Soc.}, 22(3):669--716, 2020.

\bibitem{MW}
Jason Miller and Wendelin Werner.
\newblock Connection probabilities for conformal loop ensembles.
\newblock {\em Comm. Math. Phys.}, 362(2):415--453, 2018.

\bibitem{MR3602842}
Jason Miller and Hao Wu.
\newblock Intersections of {SLE} paths: the double and cut point dimension of
  {SLE}.
\newblock {\em Probab. Theory Related Fields}, 167(1-2):45--105, 2017.

\bibitem{MR2802511}
{\c{S}}erban Nacu and Wendelin Werner.
\newblock Random soups, carpets and fractal dimensions.
\newblock {\em J. Lond. Math. Soc.}, 83(3):789--809, 2011.

\bibitem{MR3901648}
Wei Qian.
\newblock Conditioning a {B}rownian loop-soup cluster on a portion of its
  boundary.
\newblock {\em Ann. Inst. Henri Poincar\'{e} Probab. Stat.}, 55(1):314--340,
  2019.
  
\bibitem{MR4221655}
Wei Qian.
\newblock Generalized disconnection exponents.
\newblock {\em Probab. Theory Related Fields}, 179(1-2):117--164, 2021.

\bibitem{MR3896865}
Wei Qian and Wendelin Werner.
\newblock The law of a point process of {B}rownian excursions in a domain is
  determined by the law of its trace.
\newblock {\em Electron. J. Probab.}, 23:Paper No. 128, 23, 2018.

\bibitem{MR3994105}
Wei Qian and Wendelin Werner.
\newblock Decomposition of {B}rownian loop-soup clusters.
\newblock {\em J. Eur. Math. Soc.}, 21(10):3225--3253, 2019.

\bibitem{MR2153402}
Steffen Rohde and Oded Schramm.
\newblock Basic properties of {SLE}.
\newblock {\em Ann. of Math.}, 161(2):883--924, 2005.

\bibitem{MR3101840}
Oded Schramm and Scott Sheffield.
\newblock A contour line of the continuum {G}aussian free field.
\newblock {\em Probab. Theory Related Fields}, 157(1-2):47--80, 2013.

\bibitem{MR2491617}
Oded Schramm, Scott Sheffield, and David~B. Wilson.
\newblock Conformal radii for conformal loop ensembles.
\newblock {\em Comm. Math. Phys.}, 288(1):43--53, 2009.

\bibitem{MR2494457}
Scott Sheffield.
\newblock Exploration trees and conformal loop ensembles.
\newblock {\em Duke Math. J.}, 147(1):79--129, 2009.

\bibitem{MR2979861}
Scott Sheffield and Wendelin Werner.
\newblock Conformal loop ensembles: the {M}arkovian characterization and the
  loop-soup construction.
\newblock {\em Ann. of Math. (2)}, 176(3):1827--1917, 2012.

\bibitem{MR0489552}
Barry Simon.
\newblock {\em The {$P(\phi )\sb{2}$} {E}uclidean (quantum) field theory}.
\newblock Princeton Series in Physics. Princeton University Press, 1974.

\bibitem{Symanzik}
Kurt Symanzik.
\newblock Euclidean quantum field theory.
\newblock In {\em Local quantum theory, R. Jost (Ed.)}. Academic Press, 1969.

\bibitem{MR3167123}
Alain-Sol Sznitman.
\newblock On scaling limits and {B}rownian interlacements.
\newblock {\em Bull. Braz. Math. Soc. (N.S.)}, 44(4):555--592, 2013.

\bibitem{Varadhan}
S.R.S. Varadhan.
\newblock Appendix to {E}uclidean quantum field theory by {K}. {S}ymanzik.
\newblock In {\em Local quantum theory, R. Jost (Ed.)}. Academic Press, 1969.

\bibitem{MR2178043}
Wendelin Werner.
\newblock Conformal restriction and related questions.
\newblock {\em Probab. Surv.}, 2:145--190, 2005.

\bibitem{MR3618142}
Wendelin Werner.
\newblock On the spatial {M}arkov property of soups of unoriented and oriented
  loops.
\newblock In {\em S\'eminaire de {P}robabilit\'es {XLVIII}}, volume 2168 of
  {\em Lecture Notes in Math.}, pages 481--503. Springer,  2016.

\bibitem{MR4466634}
Wendelin Werner and Ellen Powell.
\newblock {\em Lecture notes on the {G}aussian free field}, {\em
  Cours Sp\'{e}cialis\'{e}s}, 28,
\newblock Soci\'{e}t\'{e} Math\'{e}matique de France, 2021.

\bibitem{MR3035764}
Wendelin Werner and Hao Wu.
\newblock From {${\rm CLE}(\kappa)$} to {${\rm SLE}(\kappa,\rho)$}'s.
\newblock {\em Electron. J. Probab.}, 18(36):1--20, 2013.

\end{thebibliography}

\end {document}